\documentclass[10pt,a4paper]{article}
\usepackage[top=3.0cm,bottom=3.0cm,left=2.5cm,right=2.5cm]{geometry}
\def\Dj{\hbox{D\kern-.73em\raise.30ex\hbox{-}
\raise-.30ex\hbox{}}}
\def\dj{\hbox{d\kern-.33em\raise.80ex\hbox{-}
\raise-.80ex\hbox{\kern-.40em}}}
\usepackage{epsfig,enumerate}
\usepackage{amsmath,amsthm,amsfonts,amssymb,amscd,cite,graphicx,array}

\newtheorem{thm}{Theorem}[section]

\newtheorem{cor}[thm]{Corollary}
\newtheorem{defi}[thm]{Definition}
\newtheorem{lem}[thm]{Lemma}

\newtheorem{problem}[thm]{Problem}
\newtheorem{exam}[thm]{Example}

 \theoremstyle{definition}

\theoremstyle{remark}

\allowdisplaybreaks

\begin{document}

\baselineskip=0.30in

\begin{center}
{\Large \bf A survey on the skew energy of oriented
graphs\footnote{Supported by NSFC No.11371205, the ``973" program
No.2013CB834204, and PCSIRT.}}

\vspace{10mm}

{\large \bf Xueliang Li}$^1$, {\large \bf Huishu Lian}$^2$

\vspace{9mm}

\baselineskip=0.20in

$^1${\it Center for Combinatorics and LPMC-TJKLC, \\
Nankai University, Tianjin 300071, P.R. China} \\
{\tt lxl@nankai.edu.cn} \\[2mm]
$^2${\it College of Science, \\
China University of Mining and Technology, Xuzhou 221116, P.R. China} \\
{\tt lhs6803@126.com}

\vspace{6mm}


\end{center}

\vspace{6mm}

\baselineskip=0.20in

\noindent {\bf Abstract }

\vspace{3mm}

{\small Let $G$ be a simple undirected graph with adjacency matrix
$A(G)$. The energy of $G$ is defined as the sum of absolute values
of all eigenvalues of $A(G)$, which was introduced by Gutman in
1970s. Since graph energy has important chemical applications, it
causes great concern and has many generalizations. The skew energy
and skew energy-like are the generalizations in oriented graphs. Let
$G^\sigma$ be an oriented graph of $G$ with skew adjacency matrix
$S(G^\sigma)$. The skew energy of $G^\sigma$, denoted by
$\mathcal{E}_S(G^\sigma)$, is defined as the sum of the norms of all
eigenvalues of $S(G^\sigma)$, which was introduced by Adiga,
Balakrishnan and So in 2010. In this paper, we summarize main
results on the skew energy of oriented graphs. Some open problems
are proposed for further study. Besides, results on the skew
energy-like: the skew Laplacian energy and skew Randi\'{c} energy
are also surveyed at the end.}

\vspace{5mm}

\baselineskip=0.30in

\section{Introduction}

In this introductory section, we will present the related background
for introducing the concept of skew energy for oriented graphs. Some
basic definitions are also given.

Let $G$ be a simple undirected graph on order $n$ with vertex set
$V(G)$ and edge set $E(G)$. Suppose $V(G)=\{v_1,v_2,\ldots,v_n\}$.
Then the adjacency matrix of $G$ is the $n\times n$ symmetric matrix
$A(G)=[a_{ij}]$, where $a_{ij}=1$ if the vertices $v_i$ and $v_j$
are adjacent, and $a_{ij}=0$ otherwise. Let $G^\sigma$ be an
oriented graph of $G$ with an orientation $\sigma$, which assigns to
each edge of $G$ a direction so that the resultant graph $G^\sigma$
becomes an oriented graph or a directed graph. Then $G$ is called
the underlying graph of $G^\sigma$. The skew adjacency matrix of
$G^\sigma$ is the $n\times n$ matrix $S(G^\sigma)=[s_{ij}]$, where
$s_{ij}=1$ and $s_{ji}=-1$ if $\langle v_i,v_j\rangle$ is an arc of
$G^\sigma$, otherwise $s_{ij}=s_{ji}=0$. It is easy to see that
$S(G^\sigma)$ is a skew symmetric matrix.

Actually, this skew adjacency matrix of an oriented graph $G^\sigma$
was first introduced by Tutte \cite{Tutte} in 1947, where he defined
a matrix $S(G^\sigma, x)=[s_{ij}]$ with $s_{ij}=x_e$ and
$s_{ji}=-x_e$ if $e=\langle v_i,v_j\rangle$ is an arc of $G^\sigma$,
otherwise $s_{ij}=s_{ji}=0$, where $x_e$ is a formal variable. Using
this matrix, he showed that the graph $G$ has a perfect matching if
and only if $\det(S(G^\sigma, x))\not\equiv 0$. It is clear that
when we take every $x_e=1$, the resultant matrix is just
the skew adjacency matrix $S(G^\sigma)$. It was proved
(see \cite{God, Lov} for examples) that the number of perfect
matchings in a graph $G$ with $m$ edges is enumerated by
$$\frac 1 {2^m} \sum_{G^\sigma\in Ort(G)}\det(S(G^\sigma)),$$
where $Ort(G)$ denotes the set of all oriented graphs of $G$. In
most cases, it is very complicated to count the number of perfect
matchings of a graph using this formula, because $\sigma$ runs over
all the $2^m$ orientations. But it is very different if a graph has
a Pfaffian orientation. A Pfaffian orientation of $G$ is such an
orientation for the edges of $G$ under which every even cycle $C$ of
$G$ such that $G\setminus V(C)$ has a perfect matching has the
property that there are odd number of edges directed in either
direction of the cycle $C$. In other words, a Pfaffian orientation
of $G$ is such an orientation for the edges of $G$ under which every
even cycle $C$ of $G$ such that $G\setminus V(C)$ has a perfect
matching is oddly oriented. In 1961, the physicists Fisher
\cite{Fisher}, Kasteleyn \cite{Kas} and Temperley \cite{Tem} used
the Pfaffian orientations of a graph to enumerate the number of
perfect matchings in a graph. It was showed that the number of
perfect matchings in $G$ is equal to the square root of
$\det(S(G^\sigma)$, where $\sigma$ is a Pfaffian orientation of $G$.
For more on the Pfaffian orientations, we refer to Robertson,
Seymour and Thomas \cite{RST}, and Thomas \cite{Thom} (45min invited
speech at the ICM 2006).

For an oriented graph $G^\sigma$ of an undirected graph $G$, there
is the adjacency matrix $A(G^\sigma)=[a_{ij}^\sigma]$, where
$a_{ij}^\sigma=1$ if $\langle v_i,v_j\rangle$ is an arc of
$G^\sigma$, and $a_{ij}^\sigma=0$ otherwise; see \cite{CDG}. It is
clear that for the adjacency matrix of $G$, we have
$A(G)=A(G^\sigma)+A^T(G^\sigma)$; whereas for the skew adjacency
matrix of $G^\sigma$, we have
$S(G^\sigma)=A(G^\sigma)-A^T(G^\sigma)$. The adjacency matrix $A(G)$
of a graph $G$ has been well-studied. It has a characteristic
polynomial, and it is symmetric and therefore all its eigenvalues
are real. Among many topics about it, we mention the energy of a
graph $G$, which was introduced by Gutman in\cite{G} and is defined
as the sum of the absolute values of all eigenvalues of $A(G)$.

Since in theoretical chemistry, the energy of a given molecular
graph is related to the total $\pi$-electron energy of the molecule
represented by that graph. Consequently, graph energy has some
specific chemistry interests and has been extensively studied. We
refer to the survey \cite{GLZ} and the book \cite{LSG} for details.
Up to now, there are various generalizations of graph energy, which
can be divided into two classes. One class contains some energy-like
quantities in undirected graphs, such as Laplacian energy, signless
Laplacian energy, incidence energy, distance energy and so on. The
other class contains skew energy and skew energy-like in oriented
graphs, such as skew Randi\'{c} energy, skew Laplacian energy and so
on. This generalization of graph energy to oriented graphs is
natural, since there are situations when chemists use digraphs
rather than graphs. One such situation is when vertices represent
distinct chemical species and arcs represent the directions in which
a particular reaction takes place between the two corresponding
species. So, people hope that skew energy and skew energy-like will
have similar applications as energy in chemistry. This paper is to
survey the main results on skew energy of oriented graphs and also
contains the results on skew Laplacian energy and skew Randi\'{c}
energy at the end.

First of all, we recall some definitions. The characteristic polynomial
of skew adjacency matrix $S(G^\sigma)$, i.e. $\det(\lambda I_n-S(G^\sigma))$,
is said to be the skew characteristic polynomial of the oriented graph
$G^\sigma$, denoted by $\phi_s(G^\sigma,\lambda)$. From linear
algebra it is known that the eigenvalues of $S(G^\sigma)$ are just
the solutions of the equation $\phi_s(G^\sigma,\lambda)=0$, which
form the spectrum of $S(G^\sigma)$ and are said to be the
skew spectrum of $G^\sigma$, denoted by $Sp_s(G^\sigma)$. Since $S(G^\sigma)$
is skew symmetric, every eigenvalue of $S(G^\sigma)$ is a pure imaginary number or 0.
The skew spectral radius of $G^\sigma$, denoted by $\rho_s(G^\sigma)$, is defined
to be the spectral radius of $S(G^\sigma)$, i.e. the maximum norm of its all eigenvalues.

Analogous to the definition of the energy of a simple undirected
graph, Adiga, Balakrishnan and So \cite{ABC} gave the definition of
skew energy of oriented graphs as follows. The definitions of the
skew Randi\'{c} energy \cite{GHL} and skew Laplacian energy
\cite{AK,AS,CLS} will be given in the corresponding sections.
\begin{defi}
Let $G^\sigma$ be an oriented graph with the skew adjacency matrix $S(G^\sigma)$
and skew spectrum $\{i\lambda_1,i\lambda_2,\ldots,i\lambda_n\}$. Then the skew energy
of $G^\sigma$, denoted by $\mathcal{E}_S(G^\sigma)$, is defined as
$$\mathcal{E}_S(G^\sigma)=\sum_{i=1}^n|\lambda_i|.$$
\end{defi}

The rest of this chapter is organized as follows: the skew energy
spans from Section \ref{matrix} to Section \ref{random} and the
final two sections are devoted to skew Randic energy and skew
Laplacian energy, respectively. Particularly, in Section
\ref{matrix}, we summarize the results about the skew characteristic
polynomial of oriented graphs. Then in Section \ref{properties}, we
collect some basic properties of the skew energy of oriented graphs
and state the integral formulas for the skew energy. Section
\ref{skew-equienergetic} is used to construct non-cospectral skew
equienergetic oriented graphs for all order $n\geq6$ and Section
\ref{bipartite} is used to survey the results about oriented graphs
with $Sp_s(G^\sigma)=iSp(G)$. General bounds of skew energy for
oriented graphs are given in Section \ref{bound}, in which the
progress on characterizing the oriented graphs achieving the upper
bound are also described. The skew spectra of various products of
graphs are determined in Section \ref{product} and therein their
applications in skew energy are illustrated. Section \ref{extremal}
outlines some results on extremal oriented graphs with regard to
skew energy. Section \ref{random} is concerned with the skew energy
of random oriented graphs. The results on the skew Randi\'{c} energy
and skew Laplacian energy of oriented graphs are listed in Section
\ref{skew-Randic} and Section \ref{skew-laplacian}, respectively. In
some sections, we propose open problems for further study.

\section{Skew characteristic polynomials of oriented graphs}\label{matrix}
This section is used to summarize the results about the
skew characteristic polynomials and skew spectra of oriented graphs.

First we give some definitions. Let $G$ be an undirected graph.
An $r$-matching of $G$ is an edge subset of $r$ edges such that
every vertex of $G$ is incident with at most one edge in it. Denote by
$m(G,r)$ the number of all $r$-matchings of $G$. Let $G^\sigma$ be
an oriented graph of $G$ and $C$ be an undirected even cycle of $G$. Then $C$
is said to be {\it evenly oriented} relative to $G^\sigma$ if it has
an even number of edges oriented in clockwise direction (and now it
also has an even number of edges oriented in anticlockwise
direction, since $C$ is an even cycle); otherwise $C$ is {\it oddly
oriented}.

Recall that a linear subgraph $L$ of $G$ is a disjoint union of some
edges and some cycles in $G$. A linear subgraph $L$ of $G$ is called
{\it even linear} if $L$ contains no odd cycle, i.e., the number of
vertices of $L$ is even. Denote by $\mathcal{EL}_i(G)$ the set of
all evenly linear subgraphs of $G$ with $i$ vertices. These
definitions in undirected graphs are the same as that in oriented
graphs. For a linear subgraph $L\in \mathcal{EL}_i$, we denote by
$p_e(L)$ and $p_o(L)$ the number of evenly oriented cycles and oddly
oriented cycles in $L$ relative to $G^\sigma$, respectively.

Let $G^\sigma$ be an oriented graph of a graph $G$ with the
skew adjacency matrix $S(G^\sigma)$. Let the skew characteristic
polynomial $\phi_s(G^\sigma,\lambda)$ of $G^\sigma$ be
\begin{equation*}
\phi_s(G^\sigma,\lambda)=\det(\lambda I_n-S(G^\sigma))=\sum^n_{i=0}c_i\lambda^{n-i}.
\end{equation*}
The following theorem, obtained by Hou and Lei in \cite{HL},
characterizes the coefficients of skew characteristic polynomial of
an oriented graph, which is analogous to the famous Sachs Theorem
\cite{Sach} for an undirected graph.
\begin{thm}\cite{HL}\label{coefficient}
Let $G^\sigma$ be an oriented graph of a graph $G$ with
the skew characteristic polynomial $\phi_s(G^\sigma,\lambda)=\sum^n_{i=0}c_i\lambda^{n-i}$. Then
\begin{equation*}
c_i=\sum_{L\in \mathcal{EL}_i}(-2)^{p_e(L)}2^{p_o(L)},
\end{equation*}
where $p_e(L)$ and $p_o(L)$ are the number of evenly oriented
cycles and the number of oddly oriented cycles of $L$ relative
to $G^\sigma$, respectively. In particular, $c_0=1$ and $c_i=0$ for all odd $i$.
\end{thm}

Gong, Li and Xu \cite{GLX} proved that all the coefficients of
skew characteristic polynomials are nonnegative.
\begin{thm}\cite{GLX}\label{coefficient2}
Let $G^\sigma$ be an oriented graph of a graph $G$ with the
skew characteristic polynomial $\phi_s(G^\sigma,\lambda)=\sum_{i\geq
0}c_{2i}\lambda^{n-2i}$. Then $c_{2i}\geq 0$ for every $i$ with
$0\leq i\leq \lfloor \frac n 2 \rfloor$.
\end{thm}

The following results are obtained by Hou and Lei in \cite{HL} as
applications of Theorem \ref{coefficient}, which can be used to find
recursions for the characteristic polynomials of some skew adjacency
matrices.

\begin{cor}\cite{HL}
Let $e=\langle u,v\rangle$ be an arc of $G^\sigma$. Then
\begin{equation*}
\begin{split}
\phi_s(G^\sigma,\lambda)=&\phi_s(G^\sigma-e,\lambda)+\phi_s(G^\sigma-u-v,\lambda)+2\sum_{e\in C\in Od(G^\sigma)}\phi_s(G^\sigma-C,\lambda)\\
&-2\sum_{e\in C\in Ev(G^\sigma)}\phi_s(G^\sigma-C,\lambda),
\end{split}
\end{equation*}
where $Od(G^\sigma)$ and $Ev(G^\sigma)$ denote the set of all oddly
oriented cycles and evenly oriented cycles of $G^\sigma$,
respectively.
\end{cor}
\begin{cor}\cite{HL}
Let $e=\langle u,v\rangle$ be an arc of $G^\sigma$ that
is on no even cycle in $G^\sigma$. Then
\begin{equation*}
\phi_s(G^\sigma,\lambda)=\phi_s(G^\sigma-e,\lambda)+\phi_s(G^\sigma-u-v,\lambda).
\end{equation*}
\end{cor}

Similarly, Xu \cite{X} got the following recursions of
the skew characteristic polynomials by deleting a vertex.
\begin{cor}\cite{X}
Let $v$ be a vertex of $G^\sigma$. Then
\begin{equation*}
\begin{split}
\phi_s(G^\sigma,\lambda)=&\lambda\phi_s(G^\sigma-v,\lambda)+\sum_{uv\in G}\phi_s(G^\sigma-u-v,\lambda)+2\sum_{v\in C\in Od(G^\sigma)}\phi_s(G^\sigma-C,\lambda)\\
&-2\sum_{v\in C\in Ev(G^\sigma)}\phi_s(G^\sigma-C,\lambda),
\end{split}
\end{equation*}
where $Od(G^\sigma)$ and $Ev(G^\sigma)$ denote the set of all oddly
oriented cycles and evenly oriented cycles of $G^\sigma$,
respectively.
\end{cor}
\begin{cor}\cite{X}
Let $v$ be a vertex of $G^\sigma$ that is on no even cycle in $G^\sigma$. Then
\begin{equation*}
\phi_s(G^\sigma,\lambda)=\lambda\phi_s(G^\sigma-v,\lambda)+\sum_{uv\in G}\phi_s(G^\sigma-u-v,\lambda).
\end{equation*}
\end{cor}

From Theorem \ref{coefficient}, it is easy to find that the
direction of an odd cycle has no effect on the coefficients
of the skew characteristic polynomial. Therefore, for an oriented
graph with no even cycles, its skew characteristic polynomial has
a special form as follows.
\begin{cor}\cite{LLL}\label{noeven}
Let $G$ be an undirected graph with no even cycles and $G^\sigma$
be an oriented graph of $G$. Then the skew characteristic polynomial
of $G^\sigma$ is of the form
\begin{equation*}
\phi_s(G^\sigma,\lambda)=\sum_{i=0}^{\lfloor n/2\rfloor} m(G,i)\lambda^{n-2i}.
\end{equation*}
\end{cor}

For an undirected graph $G$ with no even cycles, the above theorem
illustrates that its all oriented graphs has the same skew characteristic
polynomial and thus are all cospectral (i.e. the same skew spectrum).
In fact, the reverse also holds.
\begin{thm}\cite{CC}
The skew adjacency matrices of a graph $G$ are all cospectral
if and only if $G$ has no even cycles.
\end{thm}

It should be noted that Gong and Xu \cite{GX2} considered weighted
oriented graphs. They interpreted the coefficients of the
characteristic polynomials and also established some analogues of
recursions, which contains the above results as special cases.
Besides, Cavers et al. \cite{CC} obtained more general results by
considering weighted diagraphs which allows loops and dicycles of
length $2$.

It is known that skew spectrum of $G^\sigma$ consists of pure imaginary
numbers or $0$'s. Suppose that $Sp_s(G^\sigma)=\{i\lambda_1,i\lambda_2,
\ldots,i\lambda_n\}$,
where $\lambda_1\geq\lambda_2\geq\cdots\geq\lambda_n$. The skew spectrum of
the oriented graph $G^\sigma$ satisfies the following proposition.
\begin{thm}\cite{CLL3}
Let $\{i\lambda_1,i\lambda_2,\ldots,i\lambda_n\}$ be the skew
spectrum of $G^\sigma$, where
$\lambda_1\geq\lambda_2\geq\cdots\geq\lambda_n$. Then (1)
$\lambda_j=-\lambda_{n+1-j}$ for all $1\leq j\leq n$; (2) when $n$
is odd, $\lambda_{(n+1)/2}=0$ and when $n$ is even,
$\lambda_{n/2}\geq 0$; and (3) $\sum_{j=1}^{n}\lambda_j^2=2m$.
\end{thm}

We mention that many results have been obtained on the skew spectra
and skew spectral radii of oriented graphs, see \cite{AB,CC,CLL3,CLL4,CXZ,
HF,Shader,WZ,X,XG}. Moreover, the skew rank
of oriented graphs, which is defined as the rank of their skew adjacency
matrices, is studied in \cite{IMA,LY}.

\section{Some basic properties and integral formulas}\label{properties}

In this section, we first state some basic properties of the skew
energy of oriented graphs, which can be directly derived from
definition together with elementary knowledge of linear algebra.
Then we state the integral formulas for the skew energy in terms of
the skew characteristic polynomial, which play a key role in
determining the extremal oriented graphs for skew energy in some
given graph classes.

It is known that for any undirected graph with $m$ edges, there are
$2^m$ different orientations. It is natural to consider whether
there exist some operations on orientations that keep the skew
energy unchanged. Adiga, Balakrishnan and So \cite{ABC} found that
if one reverses the orientations of all the arcs incident with a
particular vertex of $G^\sigma$, then the resultant oriented graph
has the same spectrum with $G^\sigma$, and thus the same skew
energy.

Later, Hou, Shen and Zhang in \cite{HSZ} extended the above result
by introducing the {\it operation of switching}. Let $W$ be a vertex
subset of an oriented graph $G^\sigma$ and
$\overline{W}=V(G^\sigma)\backslash W$. Another oriented graph
$G^{\sigma'}$ of $G$, obtained from $G^\sigma$ by reversing the
orientations of all arcs between $W$ and $\overline{W}$, is said to
be obtained from $G^\sigma$ by switching with respect to $W$. Two
oriented graphs $G^\sigma$ and $G^{\sigma'}$ are said to be {\it
switching-equivalent} if $G^{\sigma'}$ can be obtained from
$G^\sigma$ by a sequence of switchings. Then we have the following
result.
\begin{thm}\cite{HSZ}\label{spectra}
Let $G^\sigma$ and $G^{\sigma'}$ be two oriented graphs of a graph
$G$. If $G^\sigma$ and $G^{\sigma'}$ are switching-equivalent, then
$G^\sigma$ and $G^{\sigma'}$ have the same spectra and hence
$\mathcal{E}_S(G^\sigma)=\mathcal{E}_S(G^{\sigma'})$.
\end{thm}

Let $T^\sigma$ be a labeled oriented tree rooted at vertex $v$. It
is showed that, through reversing the orientations of all arcs
incident at some vertices other than $v$, one can transform
$T^\sigma$ to another oriented tree $T^{\sigma'}$ in which the
orientations of all arcs go from low labels to high labels.
Accordingly, Adiga, Balakrishnan and So in \cite{ABC} deduced
some results on the skew energy of oriented trees.

\begin{thm}\cite{ABC}
The skew energy of an oriented tree is independent of its orientation.
\end{thm}
\begin{cor}\cite{ABC} \label{tree}
The skew energy of an oriented tree is the same as the energy of its
underlying tree.
\end{cor}

It follows that all results about the energy of undirected trees can
be immediately copied to the skew energy of oriented trees. In
particular, we have the following result.
\begin{cor}
If $T_n^\sigma$ is an oriented tree on $n$ vertices, then
\begin{equation*}
\mathcal{E}_S(S_n^\sigma)\leq \mathcal{E}_S(T_n^\sigma)\leq \mathcal{E}_S(P_n^\sigma)
\end{equation*}
where $S_n^\sigma$ and $P_n^\sigma$ denote an oriented star
and an oriented path with any orientation, respectively.
Equality holds if and only if the underlying tree $T_n$
satisfies that $T_n\cong S_n$ or $T_n\cong P_n$.
\end{cor}

In what follows, we will present the explicit expression of the skew energy
for an oriented cycle $C_n^\sigma$. Fix a vertex and label the vertices
of $C_n^\sigma$ successively. Reversing the arcs incident to a
vertex if necessary, we obtain a new oriented cycle with arcs going
from low labels to high labels possibly except one arc. Denote by
$C_n^-$ the oriented cycle with the same directions for all arcs and
by $C_n^+$ the oriented cycle with same directions for all arcs
except one arc. The skew energies of $C_n^-$ and $C_n^+$ can be
directly expressed.
\begin{thm}\cite{ABC}
Let $C_n^-$ be the oriented cycle with the same directions for
all arcs and $C_n^+$ be the oriented cycle with same directions
for all arcs except one arc. Then
\begin{equation*}
\begin{split}
&\mathcal{E}_S(C_n^-)=
\begin{cases}4\cot\frac{\pi}{n} & \text{\, if \,} n\equiv 0\,(\,mod\,2),\\
2\cot\frac{\pi}{2n} & \text{\, if \,} n\equiv 1\,(\,mod\,2);
\end{cases}\\
&\mathcal{E}_S(C_n^+)=
\begin{cases}4\csc\frac{\pi}{n} & \text{\, if \,} n\equiv0\,(\,mod\,2),\\
2\cot\frac{\pi}{2n} & \text{\, if \,} n\equiv 1\,(\,mod\,2).
\end{cases}
\end{split}
\end{equation*}
\end{thm}

We then give another two properties of the skew energy.
\begin{thm}\cite{ABC}
The skew energy of an oriented graph, if it is a rational
number, must be an even positive integer.
\end{thm}
\begin{thm}\cite{ABC}
Every even positive integer $2p$ is the skew energy of an oriented star.
\end{thm}

Similar to the Coulson integral formula for the energy of undirected
graphs \cite{GP}, we next present an integral formula \cite{ABC} for
the skew energy which enables one to compute the skew energy of an
oriented graph and compare the skew energy between two oriented
graphs without actually finding out the eigenvalues.
\begin{thm}\cite{ABC}
Let $\phi_s(G^\sigma,\lambda)$ be the skew characteristic
polynomial of an oriented graph $G^\sigma$ on $n$ vertices. Then we have
\begin{equation*}
\mathcal{E}_S(G^\sigma)=\frac{1}{\pi}\int_{-\infty}^{+\infty}
\left[n+\lambda\frac{\phi'_s(G^\sigma,-\lambda)}
{\phi_s(G^\sigma,-\lambda)}\right]d\lambda,
\end{equation*}
where $\phi'_s(G^\sigma,\lambda)$ is the derivative of $\phi_s(G^\sigma,\lambda)$.
\end{thm}

\begin{thm}\cite{HSZ}
Let $G^\sigma$ be an oriented graph of a graph $G$. Suppose that $\phi_s(G^\sigma,\lambda)=\sum^n_{i=0}c_i\lambda^{n-i}$ is the
skew characteristic polynomial of $G^\sigma$. Then
\begin{equation}\label{formula}
\mathcal{E}_S(G^\sigma)=\frac{1}{\pi}\int_{-\infty}^{+\infty}
\frac{1}{\lambda^2}\ln\left(1+\sum^{\lfloor\frac{n}{2}
\rfloor}_{i=1}c_{2i}\lambda^{2i}\right)d\lambda.
\end{equation}
\end{thm}

The following result, obtained by Zhu in \cite{Z}, illustrates an
integral formula for the difference of the skew energies of two
oriented graphs of order $n$.
\begin{thm}\cite{Z}
Let $\phi_s(G_1^{\sigma_1},\lambda)$ and
$\phi_s(G_2^{\sigma_2},\lambda)$ be two skew characteristic
polynomials of two oriented graphs $G_1^{\sigma_1}$ and
$G_2^{\sigma_2}$ of the same order. Then
\begin{equation*}
\mathcal{E}_S(G_1^{\sigma_1})-\mathcal{E}_S(G_2^{\sigma_2})=
\frac{2}{\pi}\int_{0}^{+\infty}\ln
\frac{\phi_s(G_1^{\sigma_1},\lambda)}{\phi_s(G_2^{\sigma_2},\lambda)}d\lambda\,\,.
\end{equation*}
\end{thm}

Besides, there is an integral expression to compare the skew
energies of two oriented graphs whose skew characteristic
polynomials satisfy a given recurrence relation, see \cite{Z} for
details.

From Theorem \ref{coefficient2} and the integral formula
(\ref{formula}), we find that $\mathcal{E}_S(G^\sigma)$ is a
strictly monotonically increasing function of these coefficients
$c_{2i}(G^\sigma)$ $(i=1,2,\ldots,\lfloor n/2\rfloor)$ for any
oriented graph $G^\sigma$. Therefore, the method of the quasi-order
relation $\preceq$, defined by Gutman and Polansky \cite{GP} on
graph energy, can be generalized to the skew energy of oriented
graphs. To be specific, let
$$\phi_s(G_1^{\sigma_1},\lambda)=\sum^{\lfloor
n/2\rfloor}_{i=0}c_{2i} (G_1^{\sigma_1})\lambda^{n-2i}
\text{\hspace{8pt}and\hspace{8pt}} \phi_s(G_2^{\sigma_2},\lambda)
=\sum^{\lfloor
n/2\rfloor}_{i=0}c_{2i}(G_2^{\sigma_2})\lambda^{n-2i}$$ be the
skew characteristic polynomials of two oriented graphs
$G_1^{\sigma_1}$ and $G_2^{\sigma_2}$ of order $n$, respectively. If
$c_{2i}(G_1^{\sigma_1})\leq c_{2i}(G_2^{\sigma_2})$ for all $1\leq
i\leq \lfloor n/2\rfloor$, then denote $G_1^{\sigma_1}\preceq
G_2^{\sigma_2}$, which implies that
$\mathcal{E}_S(G_1^{\sigma_1})\leq \mathcal{E}_S(G_2^{\sigma_2})$;
If $G_1^{\sigma_1}\preceq G_2^{\sigma_2}$ and there exists at least
one $j$ such that $c_{2j}(G_1^{\sigma_1})< c_{2j}(G_2^{\sigma_2})$,
then denote $G_1^{\sigma_1}\prec G_2^{\sigma_2}$, which implies that
$\mathcal{E}_S(G_1^{\sigma_1})< \mathcal{E}_S(G_2^{\sigma_2})$. This
quasi-ordering provides an important method in comparing the skew
energies of two oriented graphs. Combining this quasi-ordering with
Theorem \ref{coefficient}, Cui and Hou in \cite{CH} got the following result .

\begin{thm}\cite{CH}\label{oddly}
If $G$ has an orientation $\sigma$ such that every even cycle is
oddly oriented, then $G^\sigma$ has the maximal skew energy among
all orientations of $G$.
\end{thm}

Which graphs have an orientation $\sigma$ such that all even cycles
are oddly oriented ? Fisher and Little in \cite{FL} gave a
characterization for such graphs as follows: A graph has an
orientation under which every cycle of even length is oddly oriented
if and only if the graph contains no subgraph which is, after the
contraction of at most one cycle of odd length, an even subdivision
of $K_{2,3}$. Zhang and Li \cite{ZL} further showed that for a
bipartite graph $G$, there exists an orientation of $G$ such that
all even cycles are oddly oriented if and only if $G$ contains no
even subdivision of $K_{2,3}$, and moreover, if it is so, it must be
planar.

\section{Skew equienergetic oriented graphs}\label{skew-equienergetic}

It is obvious that if two oriented graphs have the same skew
spectra, then they possess the same skew energy, and it is known
that the operation of switching keeps skew spectra unchanged. Now it
is relatively interesting to construct some families of
non-cospectral skew equienergetic oriented graphs. Recently, Ramane
et al. \cite{RNGL} considered a join of two oriented graphs, and
established expressions for the skew characteristic polynomial and
for the skew energy of the join, respectively, by which they further
constructed non-cospectral skew equienergetic digraphs on $n$
vertices for all $n\geq6$.

We begin with some definitions. Let $G^\sigma$ be an oriented graph and $u$ be
a vertex of $G^\sigma$. The indegree of $u$, denoted by $id_{G^\sigma}(u)$, in
$G^\sigma$ is the number of arcs coming to $u$, and the outdegree of $u$,
denoted by $od_{G^\sigma}(u)$, is the number of arcs going out from $u$.
Let $G_1^\sigma$ and $G_2^\sigma$ be two oriented graphs. The join of
$G_1^\sigma$ to $G_2^\sigma$, denoted by $G_1^\sigma \rightarrow G_2^\sigma$,
is an oriented graph obtained from $G_1^\sigma$ and $G_2^\sigma$ by adding
an arc from each vertex of $G_1^\sigma$ to all vertices of $G_2^\sigma$.

The following result gives the expression for the skew characteristic polynomial
of $G_1^\sigma \rightarrow G_2^\sigma$ when $G_1^\sigma$ and $G^\sigma$ satisfy some
degree-constraints.
\begin{thm}\cite{RNGL}
For $i=1,2$, let $G_i^\sigma$ be an oriented graph of order $n_i$
with $id_{G_i^\sigma}(u)=od_{G_i^\sigma}(u)$ for all $u\in
G_i^\sigma$. Then the skew characteristic polynomial of $G_1^\sigma
\rightarrow G_2^\sigma$ is
\begin{equation*}
\phi_s(G_1^\sigma \rightarrow G_2^\sigma,\lambda)=\frac{\lambda^2+n_1n_2}{\lambda^2}
\phi_s(G_1^\sigma,\lambda)\phi_s(G_2^\sigma,\lambda).
\end{equation*}
\end{thm}

The skew energy of $G_1^\sigma \rightarrow G_2^\sigma$ can be immediately expressed
as follows.
\begin{cor}\cite{RNGL}
For $i=1,2$, let $G_i^\sigma$ be an oriented graph of order $n_i$
with $id_{G_i^\sigma}(u)=od_{G_i^\sigma}(u)$ for all $u\in
G_i^\sigma$. Then
\begin{equation*}
\mathcal{E}_S(G_1^\sigma \rightarrow G_2^\sigma)=\mathcal{E}_S(G_1^\sigma)
+\mathcal{E}_S(G_2^\sigma)+2\sqrt{n_1n_2}.
\end{equation*}
\end{cor}
\begin{thm}\cite{RNGL}\label{skew-equ}
If $H_1^\sigma$ and $H_2^\sigma$ are both non-cospectral skew
equienergetic oriented graphs of order $n$ such that
$id_{H_i^\sigma}(u)=od_{H_i^\sigma}(u)$ for all $u\in
V(H_i^\sigma)$, $i=1,2$, then for any oriented graph $G^\sigma$ with
$id_{G^\sigma}(v)=od_{G^\sigma}(v)$, $v\in V(G^\sigma)$, the
oriented graphs $H_1^\sigma\rightarrow G^\sigma$ and
$H_2^\sigma\rightarrow G^\sigma$ are non-cospectral skew
equienergetic.
\end{thm}

The above theorem provides a method to construct families of
non-cospectral skew equienergetic oriented graphs. In what follows,
we will give an example to show that for all $n\geq 6$, there exist
pairs of such oriented graphs of order $n$. The details can be found
in \cite{RNGL}.

\begin{figure}[h,t,b,p]
\begin{center}
\scalebox{1}[1]{\includegraphics{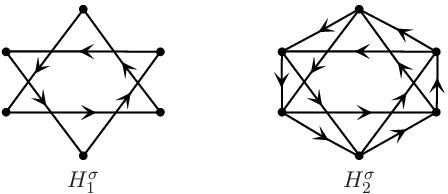}}
\end{center}
\caption{The oriented graphs $H_1^\sigma$ and $H_2^\sigma$}\label{Fig13}
\end{figure}

Take $H_1^\sigma$ and $H_2^\sigma$ as the oriented graphs depicted in Figure \ref{Fig13}.
It can be easily verified that $H_1^\sigma$ and $H_2^\sigma$ satisfy
the conditions of Theorem \ref{skew-equ} and
$$\mathcal{E}_S(H_1^\sigma)=\mathcal{E}_S(H_2^\sigma)=4\sqrt{3}.$$
Let $G^\sigma$ be any oriented graph of order $p\geq 1$ with
$id_{G^\sigma}(u)=od_{G^\sigma}(u)$, $u\in V(G^\sigma)$. Then by
Theorem \ref{skew-equ}, $H_1^\sigma\rightarrow G^\sigma$ and
$H_2^\sigma\rightarrow G^\sigma$ are non-cospectral skew
equienergetic oriented graphs. In particular,
\begin{equation*}
\mathcal{E}_S(H_1^\sigma\rightarrow G^\sigma)=
\mathcal{E}_S(H_2^\sigma\rightarrow G^\sigma)=
4\sqrt{3}+\mathcal{E}_S(G^\sigma)+2\sqrt{6p},
\end{equation*}
and $H_1^\sigma\rightarrow G^\sigma$ and $H_2^\sigma\rightarrow G^\sigma$
have the same order $n=6+p$, $p=0,1,2,\ldots$ There are many choices for $G$.
In trivial case, $G$ can be taken as $p$ isolated vertices.

The above discussions can be summarized as the following theorem.
\begin{thm}\cite{RNGL}
There exist pairs of non-cospectral skew equienergetic oriented
graphs of $n$ vertices for all $n\geq 6$.
\end{thm}

\section{Graphs with $Sp_s(G^\sigma)=iSp(G)$}\label{bipartite}

From Corollary \ref{tree} we know that for a tree $T$,
$\mathcal{E}_S(T^\sigma)=\mathcal{E}(T)$. Then, Adiga, Balakrishnan
and So \cite{ABC} posed the question that find new families of
graphs $G$ with orientations $\sigma$ such that
$\mathcal{E}_S(G^\sigma)=\mathcal{E}(G)$. Actually, for a tree $T$,
$Sp_s(T^\sigma)=iSp(T)$ for any orientation $\sigma$ of $T$. This
motivates the investigation that find new families of graphs $G$
with orientations $\sigma$ such that $Sp_s(G^\sigma)=iSp(G)$.

From Theorem \ref{coefficient}, one can deduce the following
stronger result for a tree, which was also obtained by Shader and So
in \cite{Shader} with the method of matrix analysis.
\begin{thm}\cite{HL,Shader}
Let $G$ be a graph. Then $Sp_s(G^\sigma)=iSp(G)$ for any orientation
$\sigma$ of $G$ if and only if $G$ is a tree.
\end{thm}

For bipartite graphs, we have the following results, which can also
be obtained from Theorem \ref{coefficient}.
\begin{thm}\cite{HL,Shader} \label{HL1}
A graph $G$ is bipartite if and only if there is an orientation
$\sigma$ such that $Sp_s(G^\sigma)=iSp(G)$.
\end{thm}
\begin{cor}\cite{HL} \label{HL2}
For any bipartite graph $G$, there is an orientation $G^\sigma$ with
$\mathcal{E}_S(G^\sigma)=\mathcal{E}(G)$.
\end{cor}

A natural and interesting question is that for which orientations
$\sigma$ of a bipartite graph $G$, we have $Sp_s(G^\sigma)=iSp(G)$.
It was pointed out by Shader and So in \cite{Shader} that the {\it
elementary orientation} of a bipartite graph $G=G(X,Y)$, which
assigns each edge the direction from $X$ to $Y$, is such an
orientation. Cui and Hou \cite{CH} gave a good characterization for
an oriented bipartite graph $G^\sigma$ with $Sp_s(G^\sigma)=iSp(G)$.
First recall that an even cycle $C_{2\ell}$ is said to be {\it
oriented uniformly} if $C_{2\ell}$ is oddly (resp., evenly) oriented
relative to $G^\sigma$ when $\ell$ is odd (resp., even).

\begin{thm}\cite{CH}
Let $G$ be a bipartite graph and $\sigma$ be an orientation of $G$.
Then $Sp_s(G^\sigma)=iSp(G)$ if and only if every even cycle is
oriented uniformly in $G^\sigma$.
\end{thm}

In order to ensure $Sp_s(G^\sigma)=iSp(G)$ for a given oriented
bipartite graph, the above theorem requires one to check that every
even cycle is oriented uniformly. To save the work of checking,
Chen, Li and Lian \cite{ABCLLS} proved that for an oriented bipartite
graph $G^\sigma$, $Sp_s(G^\sigma)=iSp(G)$ if and only if all chordless
cycles are oriented uniformly in $G^\sigma$. A {\it chord} of a cycle
$C$ in a graph $G$ is an edge in $E(G)\setminus E(C)$ both of whose ends
lie on $C$. A {\it chordless cycle} is a cycle without a chord. Actually,
it can be further simplified to check only a set of so-called generating
set of cycles. For details, see Remark 2.5 of \cite{ABCLLS}.

It is known from Theorem \ref{spectra} that switching-equivalence
keeps the skew spectra of an oriented graph unchanged. Therefore,
Cui and Hou \cite{CH} conjectured that such orientation $\sigma$ of
a bipartite graph $G$ that $Sp_s(G^\sigma)=iSp(G)$ is unique under
switching-equivalence. In \cite{ABCLLS} we confirmed this conjecture.

\begin{thm}\cite{ABCLLS}
Let $G=G(X,Y)$ be a bipartite graph and $\sigma$ be an orientation
of $G$. Then $Sp_s(G^\sigma)=iSp(G)$ if and only if $\sigma$ is
switching-equivalent to the elementary orientation of $G$.
\end{thm}

Some special families of oriented bipartite graphs with
$Sp_s(G^\sigma)=iSp(G)$ have been constructed \cite{CH,HL,Tian}, one
of which is obtained by considering the Cartesian product of two
oriented graphs.

Let $H$ and $G$ be graphs with $m$ and $n$ vertices, respectively.
The Cartesian product $H\square G$ of $H$ and $G$ is a graph with
vertex set $V(H)\times V(G)$ and there exists an edge between $(u_1,
v_1)$ and $(u_2, v_2)$ if and only if $u_1=u_2$ and $v_1 v_2$ is an
edge of $G$, or $v_1=v_2$ and $u_1 u_2$ is an edge of $H$. Assume
that $H^\tau$ is any orientation of $H$ and $G^\sigma$ is any
orientation of $G$.  There is a natural way to give an orientation
$H^\tau \square G^\sigma$ of $H^\tau$ and $G^\sigma$. There is an
arc from $(u_1, v_1)$ to $(u_2, v_2)$ if and only if $u_1=u_2$ and
$\langle v_1, v_2\rangle$ is an arc of $G^\sigma$, or $v_1=v_2$ and
$\langle u_1, u_2\rangle$ is an arc of $H^\tau$. It is easy to see
that $H^\tau \square G^\sigma$ is an oriented graph of $H\square G$.

\begin{thm}\cite{CH}\label{spectra2}
Let $H^\tau$ and $G^\sigma$ be oriented graphs with
$Sp_s(H^\tau)=iSp(H)$ and $Sp(G^\sigma)=iSp(G)$, respectively. Then
$Sp_s(H^\tau\square G^\sigma)=iSp(H\square G)$.
\end{thm}

Here we notice that the orientation $\sigma$ of the hypercube $Q_k$
with $Sp_s(Q_k^\sigma)=iSp(Q_k)$, obtained independently by Tian
in \cite{Tian} and Anuradha and Balakrishnan \cite{AB}, can
be viewed as a specific application of Theorem \ref{spectra2}.

It also deserves to mention that Hou and Lei in \cite{HL}
constructed from a tree the following interesting family of oriented
bipartite graphs $G^\sigma$ with $Sp_s(G^\sigma)=iSp(G)$.

Let $T$ be a tree with vertex set $\{1,2,\ldots,n\}$ and a perfect
matching $M$, i.e., $T$ is nonsingular. Let $T^\sigma$ be an oriented
tree of $T$. Suppose that $A(T)$ and $S(T^\sigma)$ are the adjacency
matrix and the skew adjacency matrix of $T$ and $T^\sigma$,
respectively. It is known that $A(T)$ and $S(T^\sigma)$ are
nonsingular since $T$ has a perfect matching. Denote by $A^{-1}(T)$
and $S^{-1}(T^\sigma)$ the inverse matrices of $A(T)$ and
$S(T^\sigma)$, respectively. A path in $T$: $P(i,j)=i_1i_2\cdots
i_{2k}$ (where $i_1=i,i_{2k}=j$ ) from a vertex $i$ to a vertex $j$
is said to be an alternating path if the edges $i_1i_2$, $i_3i_4$,
$\dots$, $i_{2k-1}i_{2k}$ are edges in the perfect matching $M$; see
\cite{BDH} for these definitions.

Define the inverse graph $T^{-1}$ of the nonsingular tree $T$ as the
graph with vertex set $\{1,2,\ldots,n\}$, where vertices $i$ and $j$
are adjacent in $T^{-1}$ if there is an alternating path between $i$
and $j$ in $T$. Suppose that $A(T^{-1})$ is the adjacency matrix of
$T^{-1}$. It is shown in \cite{BNP} that the graph $T^{-1}$ is
connected and bipartite, and the inverse matrix $A^{-1}(T)$ of the
adjacency matrix $A(T)$ of $T$ is similar to the adjacency matrix
$A(T^{-1})$ of $T^{-1}$ via a diagonal matrix of $\pm 1$.

It can be seen that the matrix
$S^{-1}(T^\sigma)$ is also skew symmetric with entries $0$ and $\pm
1$. Thus $S^{-1}(T^\sigma)$ is the skew adjacency matrix of some
oriented graph, which is defined as the inverse oriented graph of
the oriented tree $T^\sigma$, denoted by $(T^\sigma)^{-1}$. That is,
the skew adjacency matrix $S((T^\sigma)^{-1})$ of $(T^\sigma)^{-1}$
is the same as $S^{-1}(T^\sigma)$. It was shown by Hou and Lei in
\cite{HL} that $(T^\sigma)^{-1}$ is just an oriented graph of
$T^{-1}$ and the skew energy of $(T^\sigma)^{-1}$ equals the energy
of $T^{-1}$.
\begin{thm}\cite{HL}
Let $T$ be a tree with a perfect matching and $T^\sigma$ be an
oriented graph of $T$. Let $T^{-1}$ and $(T^\sigma)^{-1}$ be the
inverse graph and inverse oriented graph of $T$ and $T^\sigma$,
respectively. Then $Sp_s((T^\sigma)^{-1})=i Sp(T^{-1})$ and hence
$\mathcal{E}_S((T^\sigma)^{-1})=\mathcal{E}(T^{-1})$.
\end{thm}

Moreover, we note that Cui and Tian \cite{CT} also constructed
some families of oriented bipartite graphs with $Sp_s(G^\sigma)=iSp(G)$
by defining an operation of two oriented graphs.

\section{General bounds for the skew energy}\label{bound}

Adiga, Balakrishnan and So \cite{ABC} established a low bound and an
upper bound for the skew energy of an oriented graph $G^\sigma$ in
terms of the order and size of $G^\sigma$ as well as the maximum
degree of its underlying graph.
\begin{thm}\cite{ABC}\label{LUbound}
Let $G^\sigma$ be an oriented graph of $G$ with $n$ vertices,
m arcs and maximum degree $\Delta$. Then the skew energy
$\mathcal{E}_S(G^\sigma)$ of $G^\sigma$ satisfies that
\begin{equation*}
\sqrt{2m+n(n-1)p^{2/n}}\leq \mathcal{E}_S(G^\sigma)\leq\sqrt{2mn}\leq n\sqrt{\Delta},
\end{equation*}
where $p=|\det(S(G^\sigma))|=\prod_{i=1}^n|\lambda_i|$.
\end{thm}
\begin{cor}\cite{ABC}\label{UPB}
Any oriented graph $G^\sigma$ satisfies that $\mathcal{E}_S(G^\sigma)=n\sqrt{\Delta}$
if and only if its skew adjacency matrix satisfies that
$S(G^\sigma)^TS(G^\sigma)=\Delta I_n$, where $I_n$ is the
identity matrix of order $n$.
\end{cor}

The upper bound that $\mathcal{E}_S(G^\sigma)=n\sqrt{\Delta}$
(which, for convenience, is called the {\it optimum skew energy},
the corresponding orientation of $G$ is called the {\it optimum
orientation} and the resultant oriented graph $G^\sigma$ is called
the {\it optimum skew energy oriented graph}) implies that the
underlying graph $G$ is a $\Delta$-regular graph. Hence a natural
question was posed in \cite{ABC}:

\noindent {\bf Question:} Which $k$-regular graphs on $n$ vertices
have an orientation $G^\sigma$ with
$\mathcal{E}_S(G^\sigma)=n\sqrt{k}$, or equivalently
$S(G^\sigma)^TS(G^\sigma)=k I_n$ ?

Note that $n$ must be even since any skew symmetric matrix of odd
order has a determinant $0$, and it suffices to consider connected
$k$-regular graphs because of the following lemma.
\begin{lem}\cite{Tian}
Let $G_1^{\sigma_1}$, $G_2^{\sigma_2}$ be two disjoint oriented
graphs of order $n_1$, $n_2$ with skew adjacency matrices $S(G_1^{\sigma_1})$,
$S(G_2^{\sigma_2})$, respectively. Then for some positive integer $k$,
$S(G_1^{\sigma_1})^TS(G_1^{\sigma_1})=kI_{n_1}$ and $S(G_2^{\sigma_2})^TS(G_2^{\sigma_2})=kI_{n_2}$
if and only if the skew adjacency matrix $S(G_1^{\sigma_1} \cup G_2^{\sigma_2})$
of the union $G_1^{\sigma_1}\cup G_2^{\sigma_2}$ satisfies
$S(G_1^{\sigma_1}\cup G_2^{\sigma_2})^TS(G_1^{\sigma_1}\cup G_2^{\sigma_2})=kI_{n_1+n_2}$.
\end{lem}

Moreover, we have the following necessary condition for the complete
graph $K_n$ to have the optimum orientation.
\begin{thm}\cite{ABC}
If $K_n$ has an orientation $K_n^\sigma$ with $S(K_n^\sigma)^TS(K_n^\sigma)=(n-1)I_n$,
then $n$ is a multiple of $4$.
\end{thm}

We find that the above question is related to the so-called {\it
weighing matrices}, which are used in Combinatorial Design. A {\it
weighing matrix} $W=W(n;k)$ is defined as a square matrix with
entries $0$, $\pm 1$ having $k$ non-zero entries per row and per
column and inner product of distinct rows zero. Hence $W$ satisfies
that $WW^T=kI_n$. The number $k$ is called the weight of $W$.
Weighing matrix have been studied extensively; see
\cite{Crai,GS,MS}. Then the skew adjacency matrix of an optimum skew
energy oriented graph is a skew symmetric weighing matrices and vice
versa. Therefore, the above question is equivalent to determine the
skew symmetric weighing matrices with weighing $k$ and order $n$.

It should be pointed out that a skew symmetric weighing matrix
$W(n,n-1)$ is also called a skew symmetric conference matrix, which
is closely related to the famous Hadamard Matrix Conjecture
\cite{MS}, since $W(n,n-1)+I_n$ is a Hadamard matrix for each
skew symmetric weighing matrix $W(n,n-1)$. Moreover, there is also a
conjecture about $W(n,n-2)$ in \cite{GS}, which says that the
weighing matrices $W(4t,4t-2)$ exist for any positive integer $t$.
So far, we do not know whether the conjecture is true or not. We are
now concerned with the skew symmetric weighing matrix $W(4t,4t-2)$,
the existence of which implies the conjecture.

For small $k$, many results have been obtained. The authors in
\cite{ABC} obtained that a $1$-regular graph has an orientation with
$S(G^\sigma)^TS(G^\sigma)=I_n$ if and only if it is the graph $K_2$;
while a $2$-regular graph has an orientation with
$S(G^\sigma)^TS(G^\sigma)=2I_n$ if and only if it is the $4$-cycle
$C_4$ and has an oddly orientation.

The following lemma \cite{ABC} is very useful in characterizing
$k$-regular graphs with optimum orientations.
\begin{lem}\cite{ABC}
Let $S(G^\sigma)$ be the skew adjacency matrix of an oriented
graph $G^\sigma$. If $S(G^\sigma)^TS(G^\sigma)=kI$, then
$N(u)\cap N(v)$ is even for any two distinct vertices $u$
and $v$ of $G^\sigma$.
\end{lem}

\begin{figure}[h,t,b,p]
\begin{center}
\scalebox{1}[1]{\includegraphics{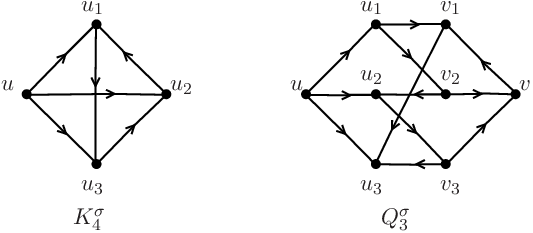}}
\end{center}
\caption{All $3$-regular oriented graphs with optimum skew energy}\label{Fig1}
\end{figure}

Applying the above lemma, Gong and Xu in \cite{GX} characterized the
underlying graphs of all $3$-regular oriented graphs with optimum
skew energy. They also gave the corresponding orientations for every
underlying graph, respectively, and further proved that such
orientation is unique for every underlying graph.
\begin{thm}\cite{GX}
Let $G^\sigma$ be a $3$-regular optimum skew energy oriented graph.
Then the underlying graph $G$ is either the complete graph $K_4$
or the hypercube $Q_3$.
\end{thm}
\begin{thm}\cite{GX}\label{3max}
Let $G^\sigma$ be a $3$-regular optimum skew energy oriented graph.
Then $G^\sigma$ (under switching-equivalent) is either $K_4^\sigma$
or $Q_3^\sigma$ depicted in Figure \ref{Fig1}.
\end{thm}

\begin{figure}[h,t,b,p]
\begin{center}
\scalebox{1}[1]{\includegraphics{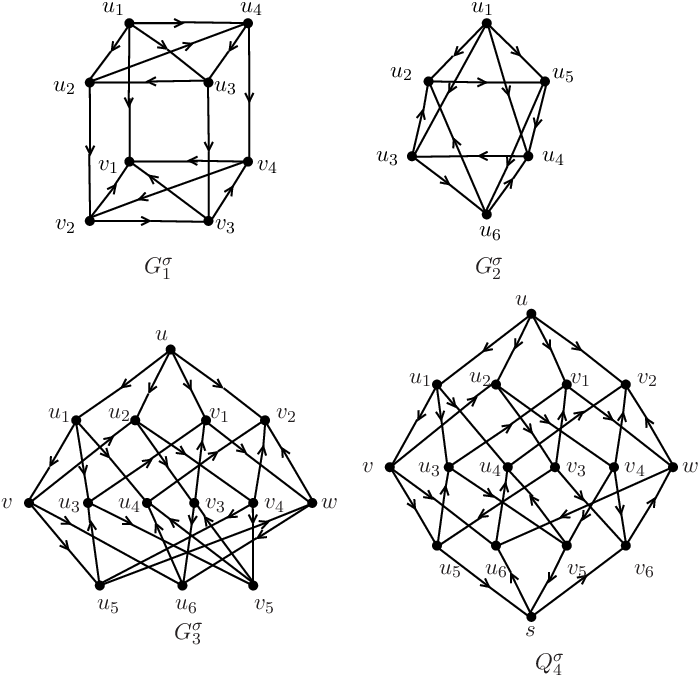}}
\end{center}
\caption{The optimum orientations for $G_1$, $G_2$, $G_3$ and $Q_4$}\label{Fig2}
\end{figure}

\begin{figure}[h,t,b,p]
\begin{center}
\scalebox{1}[1]{\includegraphics{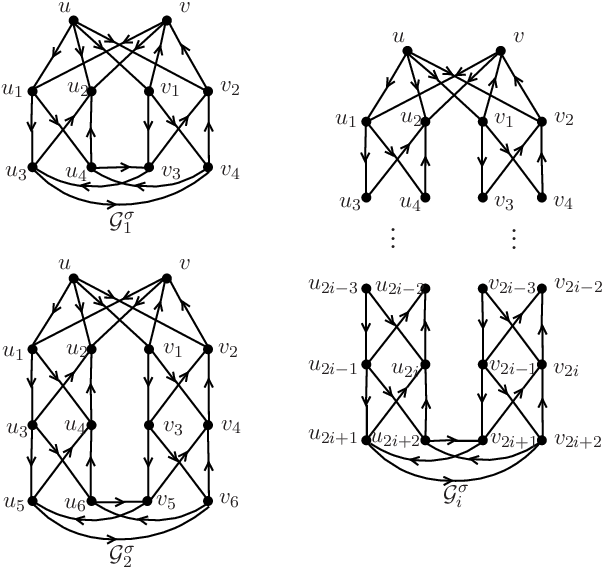}}
\end{center}
\caption{The optimum orientation for $\mathcal{G}_i$}\label{Fig3}
\end{figure}

\begin{figure}[h,t,b,p]
\begin{center}
\scalebox{1}[1]{\includegraphics{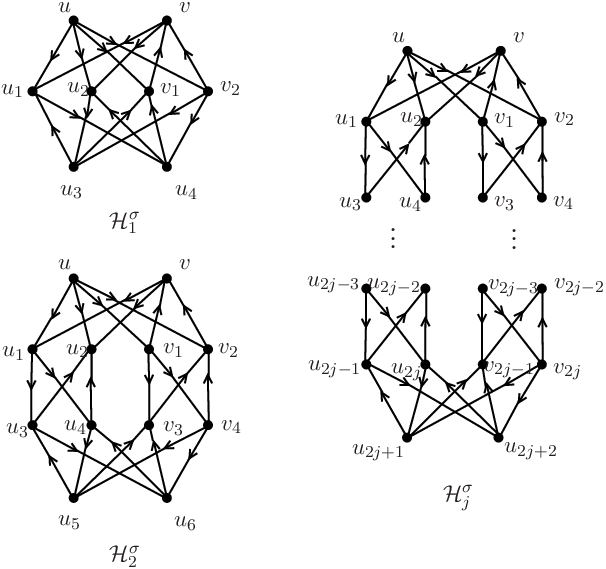}}
\end{center}
\caption{The optimum orientation for $\mathcal{H}_j$}\label{Fig4}
\end{figure}

Recently, Chen, Li and Lian in \cite{CLL2} determined the underlying
graphs of all $4$-regular oriented graphs with optimum skew energy
and gave orientations of these underlying graphs such that the skew
energies of the resultant oriented graphs indeed attain optimum.
The reader can be also referred to \cite{Lian} for details.
\begin{thm}\cite{CLL2}
Let $G$ be a $4$-regular graph. Then $G$ has an optimum orientation
if and only if $G$ is a graph of $\mathcal{F}$, which consists of
the underlying graphs of all oriented graphs in Figures \ref{Fig2},
\ref{Fig3} and \ref{Fig4}, i.e., $G_1$, $G_2$, $G_3$, the hypercube
$Q_4$, the graph $\mathcal{G}_i$ for any positive integer $i$ and
the graph $\mathcal{H}_j$ for any positive integer $j$.
\end{thm}
\begin{thm}\cite{CLL2}
Every oriented graph in Figures \ref{Fig2}, \ref{Fig3} and \ref{Fig4}
has the optimum skew energy.
\end{thm}

We point out that Gong, Zhong and Xu \cite{GZX} also independently obtained the same
result for the $4$-regular optimum skew energy oriented graph.

For any given positive integers $n$ and $k$ with $n>k\geq 5$, it is
not easy to characterize all $k$-regular graphs on order $n$ that
have an orientation $G^\sigma$ such that
$\mathcal{E}_S(G^\sigma)=n\sqrt{k}$. But it can be established that
for any positive integer $k\geq 3$, there exists a connected
$k$-regular graph $G$ that has an orientation $G^\sigma$ such that
$\mathcal{E}_S(G^\sigma)=n\sqrt{k}$, namely, the hypercube $Q_k$.
Recall that the hypercube $Q_k$ of dimension $k$ is defined
recursively in terms of the Cartesian product of graphs as follows:
\begin{equation*}
Q_k=\begin{cases} K_2, & k=1,\\ Q_{k-1}\square Q_1, & k\geq 2.
\end{cases}
\end{equation*}
Note that $Q_k$ can also be constructed by taking two copies of
$Q_{k-1}$, and then drawing an edge between each vertex in the first
copy and the corresponding vertex in the second copy. Obviously,
$Q_k$ has $2^k$ vertices and is a $k$-regular bipartite graph.
Assume that the vertex set of $Q_k$ is
$\{1,2,\ldots,2^{k-1},2^{k-1}+1,\ldots,2^k\}$. For the skew energy
of $Q_k$, Tian in \cite{Tian} obtained the following result.

\begin{thm}\cite{Tian}
For any positive integer $k$, there is an orientation of
the hypercube $Q_k$  such that the resultant oriented graph
$Q_k^\sigma$ satisfies that $\mathcal{E}_S(Q_k^\sigma)=n\sqrt{k}$.
\end{thm}

The orientation mentioned in above theorem was also given in
\cite{Tian} by the following algorithm. In fact it can be shown that
the optimum orientation for $Q_k$ is unique under
switching-equivalent.

\noindent\textbf{Algorithm:}
\begin{enumerate}[\indent Step $1.$]
\item Give the hypercube $Q_1$ an orientation $Q_1^\sigma$ such that $\langle1,2\rangle\in \Gamma(Q_1^\sigma)$.
\item Assume that $Q_1,Q_2,\dots,Q_i$
      have been oriented into $Q_1^\sigma,Q_2^\sigma,\dots,Q_i^\sigma$.
      For $Q_{i+1}$, we give an orientation $Q_{i+1}^\sigma$ using
      the following method:
\begin{enumerate}[\indent(i)]
\item Take two copies of $Q_i^\sigma$, and put an edge between
      each vertex in the first copy and the corresponding vertex in the second copy.
      Assume that the vertex set of the first copy is $\{1,2,\ldots,2^i\}$
      and the corresponding vertex set of the second copy is $\{2^i+1,2^i+2,\ldots,2^{i+1}\}$.
\item Let $\{i_1,i_2,\ldots,i_{2^{i-1}}\}$, $\{i_{2^{i-1}+1},i_{2^{i-1}+2},\ldots,i_{2^i}\}$ be a
bipartition of $\{1,2,\ldots,$ $2^i\}$. Also let
$\{i_{2^i+1},i_{2^i+2},\ldots,i_{2^i+2^{i-1}}\}$,
$\{i_{2^i+2^{i-1}+1},i_{2^i+2^{i-1}+2},\ldots,i_{2^{i+1}}\}$ be the
corresponding bipartition of $\{2^i+1,2^i+2,\ldots,2^{i+1}\}$. Give
each edge between $\{i_1,i_2,\ldots,i_{2^{i-1}}\}$ and
$\{i_{2^i+1},i_{2^i+2},\ldots,i_{2^i+2^{i-1}}\}$ an orientation such
that $\langle$ $i_k$, $i_{2^i+k}$ $\rangle$ $\in
\Gamma(Q_{i+1}^\sigma)$ for $k=1,2,\ldots,2^{i-1}$. At the same
time, give each edge between
$\{i_{2^{i-1}+1},i_{2^{i-1}+2},\ldots,i_{2^i}\}$ and
$\{i_{2^i+2^{i-1}+1},$\\$i_{2^i+2^{i-1}+2}, \ldots,i_{2^{i+1}}\}$ an
orientation such that $\langle i_{2^i+k},i_k\rangle\in
\Gamma(Q_{i+1}^\sigma)$ for $k=2^{i-1}+1,2^{i-1}+2,\ldots,2^i$.
\end{enumerate}
\item If $i+1=d$, stop; else take $i:=i+1$, return to Step 2.
\end{enumerate}

Besides the hypercube $Q_k$, many other families of oriented graphs
with optimum skew energy were characterized in \cite{CH,ABCLLS,LL},
some of which are presented in the following section.

\begin{figure}[h,t,b,p]
\begin{center}
\scalebox{0.8}[0.8]{\includegraphics{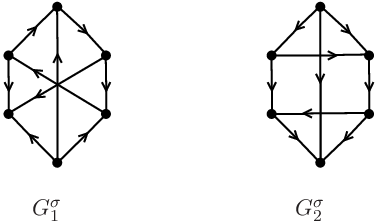}}
\end{center}
\caption{Orientations of two non-isomorphic $3$-regular graphs on
$6$ vertices having the same skew spectrum.}\label{Fig8}
\end{figure}

From the above discussions, we know that not every $k$-regular
oriented graph on $n$ vertices has an orientation such that the
resultant oriented graph attains optimum skew energy $n\sqrt{k}$.
For example, there are only two non-isomorphism $3$-regular graphs
$G_1$ and $G_2$ on $6$ vertices, both of which have no optimum
orientations. Figure \ref{Fig8} depicts the maximal orientations of
$G_1$ and $G_2$ among all orientations of $G_1$ and $G_2$. They have
the same skew spectrum $\pm 2i,\pm 2i, \pm i$ and thus the same skew
energy. But the value of their skew energy is less than $6\sqrt{3}$.
Based on this, we raise a general problem as follows:

\begin{problem}
For a given graph, how to orient it so that the oriented graph
attains the minimal or maximal skew energy ? One can try to think
about a complete graph, a hypercube, or some other graph classes.
\end{problem}

Considering that the Pfaffian orientations are important ones,
we propose the following problem:
\begin{problem}
If a graph $G$ has a Pfaffian orientation, what can we say about the
skew energy of a Pfaffian oriented graph of $G$ ?
\end{problem}

In the sequel, we consider the improvements of the bounds in Theorem \ref{LUbound}.
Note that for the lower bounds of skew energy, there are relatively few results.
Apart from the lower bound in Theorem \ref{LUbound}, Chen, Li and Lian
\cite{CLL3} obtained some better lower bounds.
\begin{thm}\cite{CLL3}
Let $G^\sigma$ be an oriented graph with $n$ vertices, $m$ arcs and
skew adjacency matrix $S$. Then
\begin{equation*}
\mathcal{E}_S(G^\sigma)\geq\sqrt{4m+n(n-2)\left(\det(S)\right)^{2/n}}.
\end{equation*}
\end{thm}

The following lower bound is derived for nonsingular oriented graphs
and the corresponding extremal graphs are also given.
\begin{thm}\cite{CLL3}
Let $G^\sigma$ be a nonsingular oriented graph with order $n$,
maximum degree $\Delta$ and skew adjacency matrix $S$. Then
\begin{equation}\label{ELB1}
\mathcal{E}_S(G^\sigma)\geq 2\sqrt{\Delta}+(n-2)\left(\frac{\det(S)}{\Delta}\right)^{\frac{1}{n-2}}
\end{equation}
and equality holds if and only if $\lambda_1=\sqrt{\Delta}$ and
$\lambda_2=\cdots=\lambda_{n/2}$.
\end{thm}

By expanding the right of Inequality (\ref{ELB1}), Chen, Li and Lian \cite{CLL3} derived
the following lower bound which is a little weaker but more simplified than the bound
(\ref{ELB1}).
\begin{cor}\cite{CLL3}
Let $G^\sigma$ be a nonsingular oriented graph with order $n$,
maximum degree $\Delta$ and skew adjacency matrix $S$. Then
\begin{equation}
\mathcal{E}_S(G^\sigma)\geq 2\sqrt{\Delta}+n-2+\ln(\det(S))-\ln{\Delta}.
\end{equation}
Equality holds if and only if $G^\sigma$ is a union of $n/2$
disjoint arcs.
\end{cor}

Recently, He and Huang \cite{HH} improved the upper bound $\sqrt{2mn}$
in Theorem \ref{LUbound} and obtained the following theorem.
\begin{thm}\cite{HH}\label{HH1}
Let $G^\sigma$ be an oriented graph with $n$ vertices, $m$ arcs and
skew adjacency matrix $S$. Then
\begin{equation*}
\mathcal{E}_S(G^\sigma)\leq\sqrt{2m(n-1)+n(\det(S))^{2/n}}.
\end{equation*}
\end{thm}

They further established a better upper bound than that in
Theorem \ref{HH1} as follows.
\begin{thm}\cite{HH}\label{HH2}
Let $G^\sigma$ be an oriented graph with $n$ vertices, $m$ arcs and
skew adjacency matrix $S$. Then
\begin{equation*}
\mathcal{E}_S(G^\sigma)\leq\sqrt{2m(n-2)+2n(\det(S))^{2/n}}.
\end{equation*}
\end{thm}

Besides, they considered the nonsingular oriented graphs, i.e., the oriented
graphs with $\det(S)\neq 0$, and derived an upper bound in terms of the order
$n$, arcs $m$ and the maximum degree $\Delta$.
\begin{thm}\cite{HH}\label{HH3}
Let $G^\sigma$ be a nonsingular oriented graph with $n$ vertices, $m$ arcs and
the maximum degree $\Delta$. Then
\begin{equation*}
\mathcal{E}_S(G^\sigma)\leq2\sqrt{\Delta}+\sqrt{(n-2)(2m-2\Delta)},
\end{equation*}
and equality holds if and only if $\lambda_1=\sqrt{\Delta}$ and
$\lambda_2=\cdots=\lambda_{n/2}$.
\end{thm}

\section{Skew-spectra and skew energies of various products of graphs}\label{product}

In this section, we consider various products of graphs, including
the Cartesian product $H\square G$, the Kronecker product $H\otimes
G$, the strong product $H\ast G$ and the lexicographic product
$H[G]$ of $H$ and $G$ where $H$ is a bipartite graph and $G$ is an
arbitrary graph. Those products are discussed in different
subsections, respectively, where we first give them orientations and
then show the skew spectra of the resultant oriented graphs. As
applications, new families of oriented graphs with optimum skew
energy are constructed in every subsection and meanwhile some
examples are given.

Let $H$ be a graph of order $m$ and $G$ be a graph of order $n$. The
definition of the Cartesian product $H\square G$ of $H$ and $G$ has
been given in Section \ref{bipartite}. Now we recall the definitions
of other graph products. The {\it Kronecker product} $H\otimes G$ of
$H$ and $G$ is a graph with vertex set $V(H)\times V(G)$ and where
$(u_1,v_1)$ and $(u_2, v_2)$ are adjacent if $u_1$ is adjacent to
$u_2$ in $H$ and $v_1$ is adjacent to $v_2$ in $G$. The {\it strong
product} $H\ast G$ of $H$ and $G$ is a graph with vertex set
$V(H)\times V(G)$; two distinct pairs $(u_1,v_1)$ and $(u_2, v_2)$
are adjacent in $H\ast G$ if $u_1$ is equal or adjacent to $u_2$,
and $v_1$ is equal or adjacent to $v_2$. The {\it lexicographic
product} $H[G]$ of $H$ and $G$ has vertex set $V(H)\times V(G)$
where $(u_1,v_1)$ is adjacent to $(u_2, v_2)$ if and only if $u_1$
is adjacent to $u_2$ in $H$, or $u_1=u_2$ and $v_1$ is adjacent to
$v_2$ in $G$.

\subsection{The orientation of $H\square G$}\label{square}

Let $H^\tau$ and $G^\sigma$ be any orientations of $H$ and $G$,
respectively. Recall the natural orientation $H^\tau\square
G^\sigma$ given in Section \ref{bipartite}. There is an arc from
$(u_1, v_1)$ to $(u_2, v_2)$ in $H^\tau \square G^\sigma$ if and
only if $u_1=u_2$ and $\langle v_1, v_2\rangle$ is an arc of
$G^\sigma$, or $v_1=v_2$ and $\langle u_1, u_2\rangle$ is an arc of
$H^\tau$. When $H$ is a bipartite graph with bipartition $X$ and
$Y$, we give another orientation of $H\square G$ by modifying the
above orientation of $H^\tau\square G^\sigma$ with the following
method. If there is an arc from $(u, v_1)$ to $(u, v_2)$ in
$H^\tau\square G^\sigma$ and $u\in Y$, then we reverse the direction
of the arc. The other arcs keep unchanged. This new orientation of
$H\square G$ is denoted by $(H^\tau\square G^\sigma)^o$. Then the
skew spectrum \cite{ABCLLS,Lian} of $(H^\tau\square G^\sigma)^o$ is given in the
following theorem.
\begin{thm}\cite{ABCLLS}\label{spectrum}
Let $H^\tau$ be an oriented bipartite graph of order $m$ and let the
skew eigenvalues of $H^\tau$ be the non-zero values $\pm \mu_1i,\pm
\mu_2i,\dots,\pm \mu_ti$ and $m-2t$ $0$'s. Let $G^\sigma$ be an
oriented graph of order $n$ and let the skew eigenvalues of
$G^\sigma$ be the non-zero values $\pm \lambda_1i,\pm\lambda_2i,\dots,
\pm \lambda_ri$ and $n-2r$ $0$'s. Then the skew eigenvalues of the
oriented graph $(H^\tau\Box G^\sigma)^o$ are $\pm i\sqrt{\mu_j^2+\lambda_k^2}$
with multiplicities $2$, $j=1,\dots,t$, $k=1,\dots,r$, $\pm\mu_ji$ with
multiplicities $n-2r$, $j=1,\dots,t$, $\pm\lambda_ki$ with multiplicities
$m-2t$, $k=1,\dots,r$, and $0$ with multiplicities $(m-2t)(n-2r)$.
\end{thm}

As an application of Theorem \ref{spectrum}, we constructed a
new family of oriented graphs with optimum skew energy in \cite{ABCLLS}.

\begin{thm} \cite{ABCLLS}\label{max-family}
Let $H^\tau$ be an oriented $\ell$-regular bipartite graph on $m$
vertices with optimum skew energy $\mathcal{E_S}(H^\tau)=m
\sqrt{\ell}$ and $G^\sigma$ be an oriented $k$-regular graph on $n$
vertices with optimum skew energy $\mathcal{E_S}(G^\sigma)=n
\sqrt{k}$. Then the oriented graph $(H^\tau\Box G^\sigma)^o$ of
$H\Box G$ has the optimum skew energy $\mathcal{E_S}((H^\tau\Box
G^\sigma)^o)=mn \sqrt{\ell+k}$.
\end{thm}

It should be noted that the special case of Theorems \ref{spectrum} and
\ref{max-family} that the bipartite graph $H$ is the
path $P_m$ of length $m$ were obtained by Cui and Hou in \cite{CH}.

It is known that there exists a $k$-regular graph with $n=2^k$ vertices
having an orientation $\sigma$ with optimum skew energy, which is the
hypercube $Q_k$. The following examples provide new families of oriented
$k$-regular graphs with optimum skew energy that have much less vertices.
\begin{exam}
Let $G_1=K_{4,4}$, $G_2=K_{4,4}\Box G_1, \dots, G_r=K_{4,4}\Box
G_{r-1}$. Because there is an orientation of $K_{4,4}$ with optimum
skew energy 16; see Figure \ref{Fig4}. Thus, we can get an orientation of
$G_r$ with optimum skew energy $2^{3r}\sqrt{4r}$. This provides a family
of $4r$-regular graphs of order $n=2^{3r}$ having an orientation
with optimum skew energy $2^{3r}\sqrt{4r}$ for $r\geq 1$.
\end{exam}
\begin{exam}
Let $G_1=K_4$, $G_2=K_{4,4}\Box G_1,\dots,G_r=K_{4,4}\Box G_{r-1}$.
It is known that $K_4$ has an orientation with optimum skew energy;
see Figure \ref{Fig1}. Thus we can get an orientation of $G_r$ with
optimum skew energy $2^{3r-1}\sqrt{4r-1}$. This provides a family of
$4r-1$-regular graphs of order $n=2^{3r-1}$ having an orientation
with optimum skew energy $2^{3r-1}\sqrt{4r-1}$ for $r\geq 1$.
\end{exam}
\begin{exam}
Let $G_1=C_4$, $G_2=K_{4,4}\Box G_1, \dots, G_r=K_{4,4}\Box
G_{r-1}$. Thus, we can get an orientation of $G_r$ with optimum
skew energy $2^{3r-1}\sqrt{4r-2}$. This provides a family of
$4r-2$-regular graphs of order $n=2^{3r-1}$ having an orientation
with optimum skew energy $2^{3r-1}\sqrt{4r-2}$ for $r\geq 1$.
\end{exam}
\begin{exam}
Let $G_1=P_2$, $G_2=K_{4,4}\Box G_1, \dots, G_r=K_{4,4}\Box
G_{r-1}$. Thus, we can get an orientation of $G_r$ with optimum
skew energy $2^{3r-2}\sqrt{4r-3}$. This provides a family of
$4r-3$-regular graphs of order $n=2^{3r-2}$ having an orientation
with optimum skew energy $2^{3r-2}\sqrt{4r-3}$ for $r\geq 1$.
\end{exam}

Moreover, if we assume that the bipartite graph $H$ is a path on two
vertices, and the graph $G$ is a tree in Theorem \ref{spectrum},
then the maximal result was obtained by Cui and Hou in \cite{CH},
whose proof used the following lemma obtained in \cite{ZY}.
\begin{lem}\cite{ZY}
Let $T$ be a tree and $\sigma$ be an arbitrary orientation of $T$.
Then the oriented graph $(P_2\square T^\sigma)^o$ has every even
cycles oddly oriented.
\end{lem}

Combining the above lemma with Theorem \ref{oddly}, the following
result is immediate.
\begin{thm}\cite{CH}
Let $T$ be a tree. Then the oriented graph $(P_2\square T^\sigma)^o$
has the maximal skew energy among all orientations of $P_2\square T$.
\end{thm}

A natural problem was proposed in \cite{CH}: for a general graph
$G$, which orientations of $P_2\square G$ yield the maximal skew
energy (or minimal skew energy) ? We are familiar with two
orientations $P_2\square G^\sigma$ and $(P_2\square G^\sigma)^o$ for
any orientation $G^\sigma$ of $G$. The following example illustrates
that $P_2\square G^\sigma$ may not be minimal among all orientations
of $P_2\square G$ when $G^\sigma$ is the minimal orientation of $G$.
But they conjectured that $(P_2\square G^\sigma)^o$ has the maximal
skew energy among all orientations of $P_2\square G$ when $G^\sigma$
is the maximal orientation of $G$.

\begin{figure}[h,t,b,p]
\begin{center}
\scalebox{1}[1]{\includegraphics{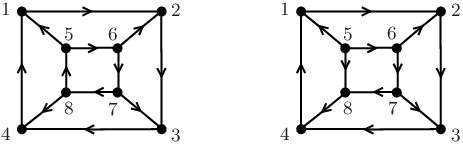}}
\end{center}
\caption{Two orientations of $P_2\square C_4$}\label{Fig9}
\end{figure}

\begin{exam}\cite{CH}
Let $G=P_2\square C_4$ (in fact, $G$ is the hypercube $Q_3$) and
$C_4^e$ be an evenly oriented. Then $Sp_s(C_4^e)=\{2i,-2i,0,0\}$ and
$\mathcal{E}_S(C_4^e)=4$ has the minimal skew energy among all
orientations of $C_4$. The left orientation in Figure \ref{Fig9},
which is $P_2\square C_4^e$, has skew energy $12$; while the right
orientation has skew energy (about $11.5$) less than $12$.
\end{exam}

\subsection{The orientation of $H\otimes G$}\label{otimes}
Let $G^{\sigma_1}_1$, $G^{\sigma_2}_2$, $G^{\sigma_3}_3$ be the
oriented graphs of order $n_1$, $n_2$, $n_3$ with skew adjacency
matrices $S_1$, $S_2$, $S_3$, respectively. Then the Kronecker product
matrix $S_1\otimes S_2\otimes S_3$ is also skew symmetric and is in fact
the skew adjacency matrix of an oriented graph of the Kronecker product
$G_1\otimes G_2\otimes G_3$. Denote the corresponding oriented graph by
$G^{\sigma_1}_1\otimes G^{\sigma_2}_2\otimes G^{\sigma_3}_3$. Adiga,
Balakrishnan and So in \cite{ABC} obtained the following result.
\begin{thm}\cite{ABC}\label{odd-Kron}
Let $G^{\sigma_1}_1$, $G^{\sigma_2}_2$, $G^{\sigma_3}_3$ be the oriented
regular graphs of order $n_1$, $n_2$, $n_3$ with optimum skew energies
$n_1\sqrt{k_1}$, $n_2\sqrt{k_2}$, $n_3\sqrt{k_3}$, respectively. Then
the oriented graph $G^{\sigma_1}_1\otimes G^{\sigma_2}_2\otimes G^{\sigma_3}_3$
has optimum skew energy $n_1n_2n_3\sqrt{k_1k_2k_3}$.
\end{thm}

It should be noted that the above Kronecker product of oriented graphs
is naturally defined, but the product requires $3$ or an odd number
of oriented graphs. In what follows, we consider the Kronecker product
of any number of oriented graphs. We first give the orientation of the
Kronecker product $H\otimes G$ where $H$ is bipartite.

Let $H$ be a bipartite graph with bipartite $(X,Y)$ and $G$ be an arbitrary
graph. For any two adjacent vertices $(u_1,v_1)$ and $(u_2,v_2)$ of $H\otimes G$,
$u_1$ and $u_2$ must be in different parts of the bipartition of vertices of $H$
and assume that $u_1\in X$. Then there is an arc from $(u_1,v_1)$ to $(u_2,v_2)$
if $\langle u_1,u_2\rangle$ is an arc of $H$ and $\langle v_1,v_2\rangle$ is an
arc of $G$, or $\langle u_2,u_1\rangle$ is an arc of $H$ and $\langle v_2,v_1\rangle$
is an arc of $G$; otherwise there is an arc from $(u_2,v_2)$ to $(u_1,v_1)$.
Denote by $(H^\tau \otimes G^\sigma)^o$ the resultant oriented graph. For the
skew spectrum of $(H^\tau \otimes G^\sigma)^o$, we obtained the following result
in \cite{LL,Lian}.
\begin{thm}\cite{LL}\label{spectrum1}
Let $H^\tau$ be an oriented bipartite graph of order $m$ and let the
skew eigenvalues of $H^\tau$ be the non-zero values $\pm \mu_1i$, $\pm
\mu_2i$, $\dots$, $\pm \mu_ti$ and $m-2t$ $0$'s. Let $G^\sigma$ be an oriented
graph of order $n$ and let the skew eigenvalues of $G^\sigma$ be the non-zero
values $\pm \lambda_1i$, $\pm\lambda_2i$, $\dots$, $\pm \lambda_ri$ and $n-2r$ $0$'s.
Then the skew eigenvalues of the oriented graph $(H^\tau\otimes G^\sigma)^o$ are
$\pm \mu_j\lambda_ki$ with multiplicities $2$, $j=1,\dots,t$, $k=1,\dots,r$,
and $0$ with multiplicities $mn-4rt$.
\end{thm}

The above theorem can be used to yield a family of oriented graphs with optimum skew energy.

\begin{thm}\cite{LL}\label{max-fami1}
Let $H^\tau$ be an oriented $k$-regular bipartite graph of order $m$
with optimum skew energy $m\sqrt{k}$. Let $G^\sigma$ be an oriented
$\ell$-regular graph of order $n$ and the optimum skew energy
$n\sqrt{\ell}$. Then $(H^\tau\otimes G^\sigma)^o$ is an oriented
$k\ell$-regular bipartite graph and has the optimum skew energy
$\mathcal{E}_S((H^\tau\otimes G^\sigma)^o)=mn\sqrt{k\ell}$.
\end{thm}

Let $H^\tau$ be an oriented bipartite graph with optimum skew
energy. Let $G_1^{\sigma_1}$ and $G_2^{\sigma_2}$ be any two
oriented graphs with the optimum skew energies. By the above
theorem, the oriented graph $(H^\tau\otimes G_1^{\sigma_1})^o$ is
bipartite and has the optimum skew energy. Therefore, the Kronecker
product $H\otimes G_1\otimes G_2$ can be oriented as
$((H^\tau\otimes G_1^{\sigma_1})^o\otimes G_2^{\sigma_2})^o$,
abbreviated as $(H^\tau\otimes G_1^{\sigma_1}\otimes
G_2^{\sigma_2})^o$, which is also bipartite and has the optimum skew
energy. The process is valid for any number of oriented graphs. Then
the following corollary is immediately implied.
\begin{cor}\cite{LL}\label{Kron}
Let $H^\tau$ be an oriented $k$-regular bipartite graph of order $m$
with optimum skew energy $m\sqrt{k}$. Let $G_i^{\sigma_i}$ be an
oriented $\ell_i$-regular graph of order $n_i$ with optimum skew
energy $n_i\sqrt{\ell_i}$ for $i=1,2,\dots,s$ and any positive
integer $s$. Then the oriented graph $(H^\tau\otimes
G_1^{\sigma_1}\otimes \cdots\otimes G_s^{\sigma_s})^o$ has the
optimum skew energy $mn_1n_2\cdots n_s\sqrt{k\,\ell_1\ell_2\cdots
\ell_s}$.
\end{cor}

It should be noted that for any even number $s$ in Corollary
\ref{Kron}, the oriented graph obtained in Corollary \ref{Kron} is
identical to the one obtained in Theorem \ref{odd-Kron}.

\subsection{The orientation of $H\ast G$}\label{ast}
Now we consider the strong product $H\ast G$ of a bipartite graph $H$
and a graph $G$. Let $H^\tau$ be an oriented graph of $H$ and $G^\sigma$
be an oriented graph of $G$. Since the edge set of $H\ast G$ is the
disjoint-union of the edge sets of $H\square G$ and $H\otimes G$, there
is a natural orientation of $H\ast G$ if $H\square G$ and $H\otimes G$
have been given orientations.

Now we give an orientation of $H\ast G$ such that the arc set of the
resultant oriented graph is the disjoint-union of the arc sets of
$(H^\tau\square G^\sigma)^o$ and $(H^\tau\otimes G^\sigma)^o$, which
are defined in the former subsections. Denote by $(H^\tau\ast
G^\sigma)^o$ this resultant oriented graph. The skew spectrum of
$(H^\tau\ast G^\sigma)^o$ is determined in the following theorem.
\begin{thm}\cite{LL}\label{spectrum2}
Let $H^\tau$ be an oriented bipartite graph of order $m$ and let the
skew eigenvalues of $H^\tau$ be the non-zero values $\pm \mu_1i$, $\pm\mu_2i$,
$\dots$, $\pm \mu_ti$ and $m-2t$ $0$'s. Let $G^\sigma$ be an oriented graph
of order $n$ and let the skew eigenvalues of $G^\sigma$ be the non-zero values
$\pm \lambda_1i$, $\pm\lambda_2i$, $\dots$, $\pm \lambda_ri$ and $n-2r$ $0$'s.
Then the skew eigenvalues of the oriented graph $(H^\tau\ast G^\sigma)^o$ are
$\pm\, i\sqrt{(u_j^2+1)(\lambda_k^2+1)-1}$ with multiplicities $2$, $j=1,\dots,t$,
$k=1,\dots,r$, $\pm\mu_ji$ with multiplicities $n-2r$, $j=1,\dots,t$, $\pm\lambda_ki$
with multiplicities $m-2t$, $k=1,\dots,r$, and $0$ with multiplicities $(m-2t)(n-2r)$.
\end{thm}

Similarly, in \cite{LL} we constructed a new family of oriented
graphs with the optimum skew energy by applying the above theorem.
\begin{thm}\cite{LL}\label{max-fami2}
Let $H^\tau$ be an oriented $k$-regular bipartite graph of order $m$
with optimum skew energy $m\sqrt{k}$. Let $G^\sigma$ be an oriented
$\ell$-regular graph of order $n$ and optimum skew energy
$n\sqrt{\ell}$. Then $(H^\tau\ast G^\sigma)^o$ is an oriented
$(k+\ell+k\ell)$-regular graph and has the optimum skew energy
$\mathcal{E}_S((H^\tau\ast G^\sigma)^o)=mn\sqrt{k+\ell+k\ell}$.
\end{thm}

Comparing Theorems \ref{max-family}, \ref{max-fami1} and
\ref{max-fami2}, we find that the oriented graphs constructed from
these theorems have the same order $mn$ but different regularities,
which are $k+\ell$, $k\ell$ and $k+\ell+k\ell$, respectively.

\begin{exam}
Let $H=C_4$, $G_0=K_4$, $G_1=H\square G_0$, \ldots, $G_r=H\square G_{r-1}$.
Obviously, $G_r$ is a $(2r+3)$-regular graph of order $4^{r+1}$. It is known
that $H$ has the orientation with the optimum skew energy $4\sqrt{2}$ and
$G_0$ has the orientation with the maximum skew energy $4\sqrt{3}$. By
Theorem \ref{max-family}, $G_r$ has the orientation with the optimum skew
energy $4^{r+1}\sqrt{2r+3}$.
\end{exam}

\begin{exam}\label{exam1}
Let $H=C_4$, $G_0=K_4$, $G_1=H\otimes G_0$, \ldots, $G_r=H\otimes G_{r-1}$.
It is obvious that $G_r$ is a $(3\cdot2^r)$-regular graph of order $4^{r+1}$.
Then by Theorem \ref{max-fami1}, $G_r$ has the orientation with optimum
skew energy $4^{r+1}\sqrt{3\cdot2^r}$.
\end{exam}

\begin{exam}\label{exam2}
Let $H=C_4$, $G_0=K_4$, $G_1=H\otimes G_0$, $G_2=G_1\otimes G_0$, \ldots,
$G_r=G_{r-1}\otimes G_0$. Note that $H$, $G_1$, $G_2$, \dots, $G_{r-1}$
are all regular bipartite graphs and $G_r$ is a $(2\cdot 3^r)$-regular
bipartite graph of order $4^{r+1}$. Then by Theorem \ref{max-fami1},
$G_r$ has the orientation with optimum skew energy $4^{r+1}\sqrt{2\cdot 3^r}$.
\end{exam}

\begin{exam}\label{exam3}
Let $H=C_4$, $G_0=K_4$, $G_1=H\ast G_0$, \ldots, $G_r=H\ast G_{r-1}$.
Note that $G_r$ is a $(4\cdot 3^r-1)$-regular graph of order $4^{r+1}$.
Then by Theorem \ref{max-fami2}, $G_r$ has the orientation with optimum skew
energy $4^{r+1}\sqrt{4\cdot 3^r-1}$.
\end{exam}

From Examples \ref{exam1}, \ref{exam2} and \ref{exam3}, one can see
that for some positive integers $k$, there exist oriented $k$-regular
graphs with the optimum skew energy, which has order $n\leq k^2$.
It is unknown that whether for any positive integer $k$, the oriented
graph exists such that its order $n$ is less than $k^2$ and it has
an orientation with the optimum skew energy.

\subsection{The orientation of $H[G]$}\label{lexi}

In this subsection, we consider the lexicographic product $H[G]$ of
a bipartite graph $H$ and a graph $G$. All definitions and notations
are the same as above. It is easy to see that the edge set $H[G]$ is
the disjoint-union of the edge sets of $H\square G$ and $H\otimes K_n$,
where $K_n$ is a complete graph of order $n$.

Let $H^\tau$ and $G^\sigma$ be oriented graphs of $H$ and $G$ with the
skew adjacency matrices $S_1$ and $S_2$, respectively. Let $K_n^\varsigma$
be an oriented graph of $K_n$ with the skew adjacency matrix $S_3$.
Then we can obtain two oriented graphs $(H^\tau\square G^\sigma)^o$ and
$(H^\tau\otimes K_n^\varsigma )^o$. Thus it is natural to yield an orientation
of $H[G]$, denoted by $H[G]^o$, such that the arc set of $H[G]^o$ is the
disjoint-union of the arc sets  of $(H^\tau\square G^\sigma)^o$ and
$(H^\tau\otimes K_n^\varsigma )$. Let $S$ be the skew adjacency matrix
of $H[G]^o$. It follows that
\begin{equation*}
S^TS=-\left[I_m\otimes S_2^2+S_1^2\otimes I_n-S_1^2\otimes S_3^2+S_1\otimes (S_2S_3)-S_1\otimes (S_3S_2)\right]
\end{equation*}

Suppose that $H^\tau$ is an oriented $k$-regular bipartite graph of
order $m$ with optimum skew energy $m\sqrt{k}$. Then
$S_1^TS_1=kI_m$. Let $G^\sigma$ be an oriented $\ell$-regular graph
of order $n$ and optimum skew energy $n\sqrt{\ell}$. Then
$S_2^TS_2=\ell I_n$. It is obvious that $H[G]^o$ is
$(kn+\ell)$-regular. Moreover, let $K_n^\varsigma$ be an oriented
graph of $K_n$ with optimum skew energy $n\sqrt{n-1}$. Then
$S_3^TS_3=(n-1)I_n$, that is, $S_3$ is a skew symmetric Hardamard
matrix \cite{MS} of order $n$. If another condition that
$S_2S_3=S_3S_2$ holds, then
$$S^TS=-\left[I_m\otimes S_2^2+S_1^2\otimes I_n-S_1^2\otimes S_3^2\right]=(kn+\ell)I_{mn}.$$
Obviously, $H[G]^o$ has the optimum skew energy $mn\sqrt{kn+\ell}$.

The following example \cite{LL} illustrates that the oriented graph
satisfying the above conditions indeed exists.
\begin{exam}\cite{LL}
Let $H^\tau$ is an arbitrary oriented $k$-regular bipartite graph of order
$m$ with the optimum skew energy $m\sqrt{k}$. Let $C_4^\sigma$ be the
oriented graph of $C_4$ with optimum skew energy $4\sqrt{2}$ and the
skew adjacency matrix $S_2$, and $K_4^\varsigma$ be the oriented graph of
$K_4$ with optimum skew energy $4\sqrt{3}$ and the skew adjacency matrix
$S_3$, see Figure \ref{Fig1}. It can be verified that $S_2S_3=S_3S_2$.
It follows that $(H[G])^o$ is an oriented $(4k+2)$-regular graph of order
$4m$ with optimum skew energy $4m\sqrt{4k+2}$.

There are many options for $H$, such as $P_2$, $C_4$, $K_{4,4}$, the hypercube
$Q_d$ and so on, which forms a new family of oriented graphs with the optimum skew energy.
\end{exam}

\section{Extremal oriented graphs in some graph classes}\label{extremal}

One of the fundamental questions that is encountered in the
study of skew energy is which oriented graphs (from a given class)
have the maximum and minimum skew energy. In this section,
we summarize these extremal results.

We have stated the results on the skew energy of oriented trees in
Section \ref{properties}. In what follows, we focus on the oriented
unicyclic graphs, oriented bicyclic graphs, oriented tricyclic
graphs, the oriented graphs with fixed numbers of vertices and arcs
and the oriented graphs without even cycles. Most of the following
results are obtained with the method of the quasi-order relation
$\preceq$ introduced in Section \ref{properties}.

\noindent{\bf Extremal oriented unicyclic graphs}

An unicyclic graph is a connected graph with the same number of
vertices and edges. Let $G(n,\ell)$ denote the set of all connected
unicyclic graphs on $n$ vertices with a unique cycle of length
$\ell$. Denote by $P_n$ and $C_n$ the path and cycle on $n$
vertices, respectively. Let $P_n^\ell$ be the unicyclic graph
obtained by connecting a vertex of $C_\ell$ with an end-vertex of
$P_{n-\ell}$, and $S_n^\ell$ be the graph obtained by connecting
$n-\ell$ pendant vertices to a vertex of $C_\ell$. Note that if $\ell$ is even,
by the switching-equivalence, there are only two different orientations on
a unicyclic graph $G$, i.e., the unique cycle $C$ is evenly oriented
or oddly oriented regardless of the orientations of other edges not
on the cycle. Denote by $G^+$ and $G^-$ the oriented graphs with the
unique cycle oddly oriented and evenly oriented, respectively. If
$\ell$ is odd, any oriented graph of a unicyclic graph $G$ has the
same skew energy, denoted by $G^\sigma$. We present some extremal
results for the skew energy of unicyclic oriented graphs as follows,
which were obtained in \cite{HSZ} by Hou, Shen and Zhang.

\begin{thm}\cite{HSZ}
For any unicyclic graph $G$, we have $G^+\succeq G^-$.
\end{thm}
\begin{thm}\cite{HSZ}
Among all orientations of unicyclic graphs on $n$ vertices,
$(S_n^3)^\sigma$ has the minimal skew energy and $(S_n^4)^-$
has the second minimal skew energy for $n\geq 6$; both $(S_5^3)^\sigma$
and $(S_5^4)^-$ have the minimal skew energy, $(S_5^4)^+$ has
the second minimal skew energy for $n=5$; $C_4^-$ has the
minimal skew energy, $(S_4^3)^\sigma$ has the second minimal skew energy for $n=4$.
\end{thm}
\begin{thm}\cite{HSZ}\label{uni-max}
Among all orientations of unicyclic graphs, $(P_n^4)^+$ is the
unique oriented graph (under switching-equivalence) with maximal
skew energy.
\end{thm}

Denote by $U_4(a,b)$ the graph obtained by attaching two pendant
paths of lengths $a$ and $b$ to the unique pendant vertex of
$S_5^4$, and then by $U^+_4(a,b)$ an oriented graph of $U_4(a,b)$
with the cycle $C_4$ oddly oriented. Moreover, let
$\mathcal{B}_n=\{U_4(a,b)|0\leq a\leq b, a+b=n-5\}$,
$\mathcal{A}_n^+=\{U^+_4(a,b)|0\leq a\leq b, a+b=n-5,a\neq 1,3,5\}$
and $\mathcal{B}_n^+=\{U^+_4(a,b)|0\leq a\leq b, a+b=n-5\}$.
Obviously, Theorem \ref{uni-max} implies that $U^+_4(0,n-5)$ has the
maximal skew energy among all oriented unicyclic graphs. The first
$\lfloor\frac{n-9}{2}\rfloor$ largest skew energies among all
oriented unicyclic graphs of order $n$ were obtained by Zhu in \cite{Z}.
\begin{thm}\cite{Z}\label{quasi}
Let $k=\lfloor\frac{n-5}{2}\rfloor$, $t=\lfloor\frac{k}{2}\rfloor$
and $\ell=\lfloor\frac{k-1}{2}\rfloor$. Then we have the following
quasi-order relation in $\mathcal{B}_n^+$:
\begin{equation*}
\begin{split}
&U^+_4(0,n-5)\succ U^+_4(2,n-7)\succ \cdots \succ U^+_4(2t,n-5-2t)\succ U^+_4(2\ell+1,n-5-2\ell-1)\\
&\succ \cdots \succ U^+_4(7,n-12)\succ U^+_4(5,n-10)\succ U^+_4(3,n-8)\succ U^+_4(1,n-6).
\end{split}
\end{equation*}
\end{thm}
\begin{thm}\cite{Z}\label{not U}
Let $G^\sigma$ be an oriented unicyclic graph with order $n\geq 31$.
If $G^\sigma\notin \mathcal{B}_n^+$, then $\mathcal{E}_S(G^\sigma)<
\mathcal{E}_S(U^+_4(7,n-12))$.
\end{thm}

Combining Theorems \ref{quasi} with \ref{not U}, we have the following result.
\begin{thm}\cite{Z}
Let $n\geq 31$. The oriented unicyclic graphs of order $n$ with the
first $\lfloor\frac{n-9}{2}\rfloor$ largest skew energies are those
in $\mathcal{A}_n^+$.
\end{thm}

Besides, Yang, Gong and Xu \cite{YGX} established the unique oriented
unicyclic graph with minimal skew energy among all oriented
unicyclic graphs on $n$ vertices with fixed diameter. Mao and Hou
\cite{MH} studied the minimal skew energy oriented unicyclic
graphs with given number of pendant vertices and girth, which can also
be found in \cite{Mao}. Gong et al.\cite{GHWXS} considered the integral
weighted oriented unicyclic graphs and determined the minimal skew energy.
Some other results on the skew energy of oriented unicyclic graphs
have been obtained by Dong \cite{D}.

\noindent {\bf Extremal oriented bicyclic graphs}

Recall that a bicyclic graph is a connected graph with $n$ vertices
and $n+1$ edges. Denote by $S_n^{3,3}$ the graph formed by joining
$n-4$ pendant vertices to a vertex of degree three of the $K_4-e$
for $n\geq 4$, by $S_n^{4,4}$ the graph formed by joining $n-5$
pendant vertices to a vertex of degree three of the complete
bipartite graph $K_{2,3}$ for $n\geq 5$, and by $P_n^{4,4}$ the
graph formed by joining two cycles of length $4$ by a path of length
$n-7$ for $n\geq 8$. For $n=6$ and $7$, $P_6^{4,4}$ and $P_7^{4,4}$
denote two new graphs. Let $P_6^{4,4}$ be the graph consisting of
two cycles of length $4$ sharing a common edge and $P_7^{4,4}$ be
the graph consisting of two cycles of length $4$ sharing a common vertex.
Let $C_x$ and $C_y$ be two cycles in a bicyclic graph $G$ with $t$ common vertices.
Note that if $t\leq 1$, then $G$ contains exactly two cycles. If $t\geq
2$, then $G$ contains exactly three cycles. The third cycle is
denoted by $C_z$, where $z=x+y-2t+2$. Without loss of generality,
assume $x\leq y\leq z$.

For convenience, when $t\leq 1$ denote by $G^{a,b}$ the oriented
bicyclic graph on which $C_x$ is of orientation $a$ and $C_y$ is of
orientation $b$, where $a,b\in\{+,-,*\}$ and $*$ denotes any
orientation; when $t\geq 2$, denote by $G^{a,b,c}$ the oriented
bicyclic graph on which $C_x$ is of orientation $a$, $C_y$ is of
orientation $b$ and $C_z$ is of orientation $c$, where
$a,b,c\in\{+,-,*\}$. Shen, Hou and Zhang \cite{HSZ2} studied the
minimal and maximal skew energy oriented bicyclic graphs and obtained
the following results.
\begin{thm}\cite{HSZ2}
Among all oriented bicyclic graphs of order $n$,
$(S_n^{3,3})^{*,*,-}$ has the minimal skew energy for $n\geq 8$;
both $(S_7^{3,3})^{*,*,-}$ and $(S_7^{4,4})^{-,-,-}$ have the
minimal skew energy for $n=7$; $(S_n^{4,4})^{-,-,-}$ has the minimal
skew energy for $n=5,6$.
\end{thm}
\begin{thm}\cite{HSZ2}
Among all oriented bicyclic graphs of order $n\geq 8$,
$(P_n^{4,4})^{+,+}$ has the maximal skew energy; $(P_7^{4,4})^{+,+}$
have the maximal skew energy for $n=7$; $(P_6^{4,4})^{+,+,+}$ has
the maximal skew energy for $n=6$.
\end{thm}

\begin{figure}[h,t,b,p]
\begin{center}
\scalebox{0.8}[0.8]{\includegraphics{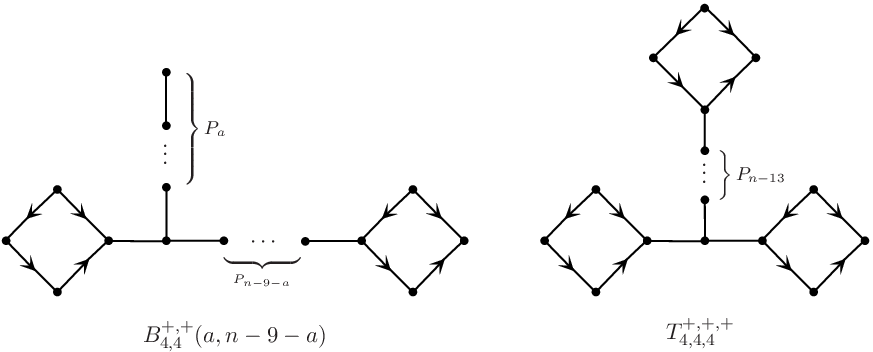}}
\end{center}
\caption{The extremal oriented graphs $B^{+,+}_{4,4}(a,n-9-a)$ and $T^{+,+,+}_{4,4,4}$}\label{Fig14}
\end{figure}

Later, Wang, Zhao and Ye \cite{WZY} established the oriented
bicyclic graphs with the second largest skew energy. Recently,
Qin, Yang and Wang \cite{QYW} further obtained the extensive result
by determining the oriented bicyclic graphs of order $n\geq13$ with
the first five largest skew energies. The extremal oriented graphs
$B_{4,4}(a,n-9-a)$ are shown in Figure \ref{Fig14},
where $0\leq a\leq n-9$.
\begin{thm}\cite{QYW}
Among all oriented bicyclic graphs with order $n\geq13$, the graphs
$B^{+,+}_{4,4}(0,n-9)\succeq B^{+,+}_{4,4}(n-9,0)\succeq
B^{+,+}_{4,4}(2,n-11)\succeq B^{+,+}_{4,4}(n-11,2)\succeq
B^{+,+}_{4,4}(4,n-13)$
have the first five largest skew energies.
\end{thm}

\noindent {\bf Extremal oriented tricyclic graphs}

A tricyclic graph is a connected graph with $n$ vertices and $n+2$ edges.
Li, Qin, Wang and Yang \cite{LQYW} characterized the oriented graph
with maximal skew energy among all oriented tricyclic graphs.
\begin{thm}
Among all oriented oriented tricyclic graphs with order $n\geq13$, the
oriented graph $T^{+,+,+}_{4,4,4}$ has the maximal skew energy, see Figure \ref{Fig14}.
\end{thm}

\noindent {\bf Extremal oriented graphs with $n$ vertices and $m$ $(n\leq m\leq 2(n-2))$ arcs}

\begin{figure}[h,t,b,p]
\begin{center}
\scalebox{1}[1]{\includegraphics{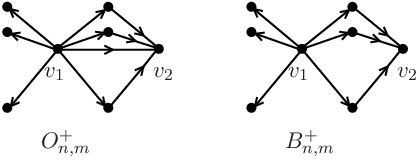}}
\end{center}
\caption{Two oriented graphs $O_{n,m}^+$ and $B_{n,m}^+$}\label{Fig7}
\end{figure}

Now we consider the oriented graphs on $n$ vertices with
$m$ $(n\leq m\leq 2(n-2))$ arcs. The oriented graphs with
minimal skew energy in this oriented graph class were
characterized in \cite{GLX} by Gong, Li and Xu. Let $O_{n,m}^+$ be
the oriented graph on $n$ vertices which is obtained from the
oriented star $S_n^\sigma$ with center $v_1$ by adding $m-n+1$ arcs
such that all those arcs have a common vertex $v_2$, where $v_1$ is
the tail of each arc incident to it and $v_2$ is the head of each
arc incident to it; and $B_{n,m}^+$ be the oriented graph obtained
from $O_{n,m+1}^+$ by deleting the arc $\langle v_1,v_2\rangle$; see
Figure \ref{Fig7}.

\begin{thm}\cite{GLX}
Let $G^\sigma$ be an oriented graph with minimal skew energy among
all oriented graphs with $n$ vertices and $m$ $(n\leq m\leq 2(n-2))$
arcs. Then, up to isomorphic, $G^\sigma$ is
\begin{itemize}
\item [(1)] $O_{n,m}^+$ if $m<\frac{3n-5}{2}$;
\item [(2)] either $B_{n,m}^+$ or $O_{n,m}^+$ if $m=\frac{3n-5}{2}$; and
\item [(3)] $B_{n,m}^+$ otherwise.
\end{itemize}
\end{thm}

\noindent {\bf Extremal oriented graphs without even cycles}

We consider the extremal skew energy of the oriented graphs
without even cycles, see \cite{Li,LLL}. Denote by $\mathcal{O}_n$ the
class of oriented graphs of order $n$ with no even cycles, and
by $\mathcal{O}_{n,m}$ the class of oriented graphs in $\mathcal{O}_n$
with $m$ arcs. From corollary \ref{noeven}, it is easy to find that
in the class $\mathcal{O}_n$ the skew characteristic polynomial of
an oriented graph is independent of its orientation, so is the skew energy.
Thus we denote by $G^*$ the oriented graph with any orientation.

In \cite{LLL}, we determined the minimal skew energy oriented graphs
in $\mathcal{O}_n$ and $\mathcal{O}_{n,m}$ with $n-1\leq
m\leq \frac{3}{2}(n-1)$. Denote by
$S_{n,m}$, $n-1\leq m\leq \frac{3}{2}(n-1)$ the graph obtained from $S_n$
by attaching $m-n+1$ edges to different pairs of pendant vertices, such that
the pendant vertices of $S_n$ is of degree no more than $2$.
Clearly, $S_{n,n-1}\cong S_n$.
\begin{thm}\cite{LLL}
Let $G^\sigma$ be an oriented graph in $\mathcal{O}_n$. If $G \ncong S_n$,
then $\mathcal{E}_s(G^\sigma)>\mathcal{E}_s(S^*_n)$.
\end{thm}
\begin{thm}\cite{LLL}
Let $G^\sigma$ be an oriented graph in $\mathcal{O}_{n,m}$, $n-1\leq
m\leq \frac{3}{2}(n-1)$. If $G\ncong S_{n,m}$, then
$\mathcal{E}_s(G^\sigma)>\mathcal{E}_s(S^*_{n,m})$.
\end{thm}

We \cite{LLL} also characterized the maximal skew energy oriented graphs in $\mathcal{O}_{n,n}$
and $\mathcal{O}_{n,n+1}$, and in the latter case we assume that $n$ is even.
\begin{thm}\cite{LLL}
If $n$ is odd, the oriented graph with maximal skew energy in
$\mathcal{O}_{n,n}$ is $(C_n)^*$. If $n\equiv 0\ (mod \ 4)$, the oriented graph
with maximal skew energy in $\mathcal{O}_{n,n}$ is $(P_n^{\ell_0})^*$,
$\ell_0=n/2+1$. If $n\equiv 2\ (mod \ 4)$, the oriented graph with maximal
skew energy in $\mathcal{O}_{n,n}$ is $(P_n^{\ell_1})^*$, $\ell_1=n/2$
or $n/2+2$.
\end{thm}
\begin{thm}\cite{LLL}
let $n$ be an even integer, if $n\equiv 0\ (mod \ 4)$, the oriented graph
with maximal skew energy in $\mathcal{O}_{n,n+1}$ is
$P_n^{\ell,\ell-2}$, $\ell=n/2+1$. If $n\equiv 2\ (mod \ 4)$, the
oriented graph with maximal skew energy in $\mathcal{O}_{n,n+1}$ is
$P_n^{\ell',\ell'}$, $\ell'=n/2$.
\end{thm}

\section{The skew energy of random oriented graphs}\label{random}

In this section, we state our results on the skew energy of random
oriented graphs \cite{CLL,Lian}. At first let us recall the random
oriented graph model and the random oriented regular graph model. A
random oriented graph $G^\sigma(n,p)$ on $n$ vertices is obtained by
drawing an edge between each pair of vertices, randomly and
independently, with probability $p$ and then orienting each existing
edge, randomly and independently, with probability $1/2$. A random
regular graph $G_{n,d}$, where $d=d(n)$ denotes the degree, is a
random graph chosen uniformly from the set of all simple $d$-regular
graphs on $n$ vertices. A random oriented regular graph, denoted by
$G_{n,d}^\sigma$, is obtained by orienting each edge of the random
regular graph $G_{n,d}$, randomly and independently, with
probability $1/2$.

Given a random graph model $\mathcal{G}(n,p)$, we say that almost
every graph $G(n,p)\in \mathcal{G}(n,p)$ has a certain property
$\mathcal{P}$ if the probability that $G(n,p)$ has the property
$\mathcal{P}$ tends to $1$ as $n\rightarrow \infty$, or we say
$G(n,p)$ {\it almost surely (a.s. )} satisfies the property
$\mathcal{P}$.

In what follows, we state the exact estimates of the skew energy for
almost all oriented graphs and almost all oriented regular graphs,
respectively, where the estimate of the latter is distinguished into
two cases according to the value of regular degree $d$; see
\cite{CLL} for details.
\begin{thm}\cite{CLL}
For $p=\omega(\frac{1}{n})$, the skew energy
$\mathcal{E}_S(G^\sigma(n,p))$ of the random oriented graph
$G^\sigma(n,p)$ enjoys a.s. the following formula:
\begin{equation*}
\mathcal{E}_S(G^\sigma(n,p))=n^{3/2}p^{1/2}\left(\frac{8}{3\pi}+o(1)\right).
\end{equation*}
\end{thm}

Since $p$ varies in the interval $[0,1]$, we can get the following
result.
\begin{cor}\label{cor1}
It is almost sure that the skew energy is increasing as the number
of arcs is getting large.
\end{cor}

This also verifies that a tournament (an oriented complete graph)
has the maximum skew energy.

It is known from \cite{DLL, LSG} that the energy of an undirected
random graph $G(n,p)$ enjoys the formula below:
\begin{equation*}
\mathcal{E}(G(n,p))=n^{3/2}[p(1-p)]^{1/2}\left(\frac{8}{3\pi}+o(1)\right).
\end{equation*}
Since $[p(1-p)]^{1/2}$ is less than $p^{1/2}$ for any $0<p<1$, the
following corollary is immediate.
\begin{cor}\label{cor2}
It is almost sure that for every graph $G$, the energy of $G$ is
less than the skew energy of an oriented graph $G^\sigma$ of $G$.
\end{cor}

\begin{thm}\cite{CLL}
For any fixed integer $d\geq 2$, the skew energy
$\mathcal{E}_S(G^\sigma_{n,d})$ of the random oriented regular graph
$G^\sigma_{n,d}$ enjoys a.s. the following formula:
\begin{equation*}
\mathcal{E}_S(G^\sigma_{n,d})=n\left(\frac{2d\sqrt{d-1}}{\pi}-\frac{d(d-2)}{\pi}\cdot \arctan\frac{2\sqrt{d-1}}{d-2}+o(1)\right).
\end{equation*}
In particular, when $d=2$,
$\mathcal{E}_S(G^\sigma_{n,d})=n\left(4/\pi+o(1)\right)$.
\end{thm}

\begin{thm}\cite{CLL}
For $d=d(n)\rightarrow \infty$, the skew energy
$\mathcal{E}_S(G^\sigma_{n,d})$ of the random oriented regular graph
$G^\sigma_{n,d}$ enjoys a.s. the following formula:
\begin{equation*}
\mathcal{E}_S(G^\sigma_{n,d})\,=\,nd^{1/2}\left(\frac{8}{3\pi}+o(1)\right).
\end{equation*}
\end{thm}

\noindent{\bf Remark 9.1.} Corollary \ref{cor1} holds only ``almost surely", but not
``always surely".

There are oriented graphs $G^\sigma$ such that some of their proper
subgraphs have larger skew energy than that of $G^\sigma$. From page
57 of \cite{LSG} we know that, as the energy of undirected graphs is
concerned, $\mathcal{E}(K_{n,n})<\mathcal{E}(K_{n,n}\setminus e)$
for any edge of the complete bipartite graph $K_{n,n}$. Then, from
Theorem \ref{HL1}, Corollary \ref{HL2} and the context thereafter,
we know that there is an orientation $\sigma$ of the bipartite graph
$K_{n,n}$ such that
$\mathcal{E_S}(K_{n,n}^\sigma)=\mathcal{E}(K_{n,n})$ and
$\mathcal{E_S}((K_{n,n}\setminus
e)^\sigma)=\mathcal{E}(K_{n,n}\setminus e)$, which gives that
$\mathcal{E_S}(K_{n,n}^\sigma)<\mathcal{E_S}((K_{n,n}\setminus
e)^\sigma)$.

\noindent{\bf Remark 9.2.} Corollary \ref{cor2} holds also only ``almost surely", but not
``always surely".

There are undirected graphs $G$ such that some of their oriented
graphs $G^\sigma$ have smaller skew energy than the energy of $G$.
From page 26 of \cite{LSG} we know that the energy of an undirected
cycle $C_n$ is as follows:
\begin{equation*}
\mathcal{E}(C_n)=
\begin{cases}4\cot\frac{\pi}{n} & \text{\, if \,} n\equiv 0\,(\,mod\,4),\\
4\csc\frac{\pi}{n} & \text{\, if \,} n\equiv 2\,(\,mod\,4),\\
2\csc\frac{\pi}{2n} & \text{\, if \,} n\equiv 1\,(\,mod\,2).
\end{cases}
\end{equation*}
Comparing the energy of the undirected cycle $C_n$ and the skew
energies of the two oriented cycles $C_n^+$ and $C_n^-$, we obtain
the following result:
\begin{equation*}
\begin{split}
\mathcal{E}(C_n)=\mathcal{E}_S(C_n^-)<\mathcal{E}_S(C_n^+)\hspace{15pt}&\text{for}\hspace{7pt}n\equiv 0\,(\,mod\,4),\\
\mathcal{E}_S(C_n^-)<\mathcal{E}_S(C_n^+)=\mathcal{E}(C_n)\hspace{15pt}&\text{for}\hspace{7pt}n\equiv 2\,(\,mod\,4),\\
\mathcal{E}_S(C_n^-)=\mathcal{E}_S(C_n^+)<\mathcal{E}(C_n)\hspace{15pt}&\text{for}\hspace{7pt}n\equiv
1\,(\,mod\,2).
\end{split}
\end{equation*}
From the above one can see that there are cases that the energy of
an undirected cycle is smaller than the skew energy of its oriented
cycle, and there are some other cases that the skew energy of an
oriented cycle is smaller than the energy of the cycle.

\section{Shew Randi\'{c} energy of oriented graphs}\label{skew-Randic}
The skew Randi\'{c} energy of oriented graphs \cite{GHL} can be viewed as
the generalization of Randi\'{c} energy of undirected graphs. This
section is devoted to collecting the results of the skew Randi\'{c}
energy of oriented graphs, most of which are analogous to those of
skew energy of oriented graphs.

Recall first the Randi\'{c} energy of undirected graphs. Let $G$ be an
undirected graph with vertex set $V(G)=\{v_1,v_2,\ldots,v_n\}$ and let
$d_i$ be the degree of $v_i$, $i=1,2,\ldots,n$. The Randi\'{c} matrix of
$G$ is the $n\times n$ matrix $R(G)=[r_{ij}]$, where $r_{ij}=(d_id_j)^{-1/2}$
if $v_i$ and $v_j$ are adjacent, and $r_{ij}=0$ otherwise. The Randi\'{c}
matrix $R(G)$ can be regarded as the adjacency matrix of a weighted graph
and its characteristic polynomial is called the $R$-characteristic polynomial
of $G$. Obviously, $R(G)$ is symmetric and thus all eigenvalues are real,
which form the Randi\'{c} spectrum $Sp_{R}(G)$ of $G$. The Randi\'{c} energy of $G$
is defined as the sum of the absolute values of all eigenvalues of $S(G)$,
denoted by $RE(G)$.

Let $G^\sigma$ be an oriented graph of $G$. The skew Randi\'{c} matrix of
$G^\sigma$ is defined as $R_s(G^\sigma)=[\hat{r}_{ij}]$, where
$\hat{r}_{ij}=(d_id_j)^{-1/2}$ and $\hat{r}_{ji}=-(d_id_j)^{-1/2}$ if
$(v_i,v_j)$ is an arc of $G^\sigma$, $\hat{r}_{ij}=\hat{r}_{ji}=0$ otherwise.
Note that $R_s(G^\sigma)$ can be viewed as the skew adjacency matrix of a
weighted oriented graph. The characteristic polynomial of $R_s(G^\sigma)$,
i.e. $\phi_{R_s}(G^\sigma,\lambda)=\det(\lambda I_n-S(G^\sigma))$, is said
to be the $R_s$-characteristic polynomial of $G^\sigma$. It is easy to
check that $R_s(G^\sigma)$ is skew symmetric and thus its spectrum consists
of purely imaginary numbers or $0$'s, which is also called skew Randi\'{c}
spectrum of $G^\sigma$ and is denoted by $Sp_{R_s}(G^\sigma)$. The skew
Randi\'{c} energy is defined as the sum of norms of all eigenvalues of
$R_s(G^\sigma)$, denoted by $RE_S(G^\sigma)$.

We first consider $R_s$-characteristic polynomial and the skew Randic spectra
of $G^\sigma$, which are similar to the skew characteristic polynomial and
skew spectra of $G^\sigma$.
\begin{thm}\cite{GX2,GHL}
Let $G^\sigma$ be an oriented graph of a graph $G$ with the $R_s$-characteristic polynomial
$\phi_{R_s}(G^\sigma,\lambda)=\sum^n_{i=0}c_i\lambda^{n-i}$. Then
\begin{equation*}
c_i=\sum_{L\in \mathcal{EL}_i}(-2)^{p_e(L)}2^{p_o(L)}W(L),
\end{equation*}
where $p_e(L)$ and $p_o(L)$ are the number of evenly oriented cycles and
the number of oddly oriented cycles of $L$ relative to $G^\sigma$, respectively,
$W(L)$ are the product of the weights of all arcs of $L$. In particular,
(i) $c_0=1$, (ii) $c_2=R_{-1}(G)=\sum_{v_i\sim v_j}(d_id_j)^{-1}$,
the general Randic index with $\alpha=-1$, (iii) $c_i=0$ for all odd $i$.
\end{thm}

\begin{thm}\cite{GHL}
Let $\{i\lambda_1,i\lambda_2,\ldots,i\lambda_n\}$ be the skew Randi\'{c}
spectrum of $G^\sigma$, where $\lambda_1\geq\lambda_2\geq\ldots\geq\lambda_n$.
Then (1) $\lambda_i=-\lambda_{n-1-i}$ for all $1\leq i\leq n$; (2)
when $n$ is odd, $\lambda_{(n+1)/2}=0$ and when $n$ is even, $\lambda_{n/2}\geq0$;
and (3) $\sum_{i=1}^n\lambda_i^2=2R_{-1}(G)$, where $R_{-1}(G)$ is the
general Randi\'{c} index of $G$ with $\alpha=-1$.
\end{thm}

Next, we list the properties of skew Randi\'{c} energy of oriented graphs,
which are same to those of skew energy of oriented graphs. The following
theorem shows that two switching-equivalent oriented graphs possess the
same skew Randi\'{c} spectra.
\begin{thm}\cite{GHL}
If $G^\sigma$ and $G^\tau$ are switching-equivalent, then $Sp_{R_s}(G^\sigma)=Sp_{R_s}(G^\tau)$
\end{thm}

By the above theorem, we can deduce the propositions of skew Randi\'{c}
energy of oriented trees.
\begin{thm}\cite{GHL}
The skew Randi\'{c} energy of an oriented tree is independent of its orientation.
\end{thm}
\begin{cor}\cite{GHL}
The skew Randi\'{c} energy of an oriented tree is the same as the Randi\'{c}
energy of its underlying tree.
\end{cor}

Moreover, we consider the oriented graphs with $Sp_{R_s}(G^\sigma)=iSp_R(G)$.
\begin{thm}\cite{GHL}
$G$ is a forest if and only if for any orientation $\sigma$ of $G$,
$Sp_{R_s}(G^\sigma)=iSp_R(G)$.
\end{thm}
\begin{thm}\cite{GHL}
$G$ is a bipartite graph if and only if there is an orientation $\sigma$ of $G$
such that $Sp_{R_s}(G^\sigma)=iSp_R(G)$.
\end{thm}
\begin{thm}\cite{GHL}
Let $G$ be a bipartite graph and $\sigma$ be an orientation of $G$. Then
$Sp_{R_s}(G^\sigma)=iSp_R(G)$ if and only if $\sigma$ is switching-equivalent
to the elementary orientation of $G$.
\end{thm}

For the remainder of this section, we consider the bounds of skew Randi\'{c}
energy of oriented graphs.
\begin{thm}\cite{GHL}
$\sqrt{4R_{-1}(G)+n(n-2)p^{\frac{2}{n}}}\leq RE_S(G^\sigma)\leq 2\sqrt{\lfloor\frac{n}{2}\rfloor R_{-1}(G)}$, where $p=|\det(R_S(G^\sigma))|$.
\end{thm}

Note that there are many results on the upper and lower bounds on $R_{-1}(G)$.
Combining with the above theorem, we can obtain the upper and lower bounds on
the skew Randi\'{c} energy without the parameter $R_{-1}(G)$, one of which is
given as follows.
\begin{thm}\cite{GHL}
Let $G^\sigma$ be an oriented graph of order $n$ without no isolated vertices. Then
$$\sqrt{\frac{2n}{n-1}+n(n-2)p^{\frac{2}{n}}}\leq RE_S(G^\sigma)\leq 2\lfloor\frac{n}{2}\rfloor,$$
where $p=|\det(R_s(G^\sigma))|$. The equality in the lower bound holds
if and only if $G$ is a complete graph with exact two nonzero skew Randi\'{c}
eigenvalues when $n$ is odd, and $R_s(G^\sigma)^TR_s(G^\sigma)=\frac{1}{n-1}I_n$
when $n$ is even. The equality in the upper bound holds if and only if either
$n$ is even and $G$ is the disjoint union of paths of length $1$, or $n$ is
odd and $G$ is the disjoint union of $\frac{n-3}{2}$ paths of length $1$ and
one path of length $2$, and $\sigma$ is an arbitrary orientation of $G$.
\end{thm}

We point out that the bounds of skew Randi\'{c} energy of oriented
trees and oriented chemical trees have been determined, and the
corresponding extremal graphs have also been characterized; see
\cite{GHL} for details.

\section{Skew Laplacian energy of oriented graphs}\label{skew-laplacian}

In previous sections, we summarize the results of skew energy of
oriented graphs. This section is concerned with the skew Laplacian
energy of oriented graphs.

Before proceeding, we will briefly introduce the Laplacian energy of
undirected graphs in order to do comparing. We mention first that
there are two different definitions of the Laplacian energy. Let $G$
be a simple undirected graph with vertex set
$\{v_1,v_2,\ldots,v_n\}$ and $m$ edges. Let $A(G)$ be the adjacency
matrix of $G$ and $D(G)=\mbox{diag}(d_1,d_2,\ldots,d_n)$ be the
diagonal matrix of vertex degrees. Then Denote by $L(G)=D(G)-A(G)$
the Laplacian matrix of $G$ with spectrum
$\{\mu_1,\mu_2,\ldots,\mu_n\}$.

Kragujevac \cite{Kra} gave a definition for Laplacian energy using the second
spectral moment, namely, $LE_k(G)=\sum_{i=1}^n\mu_i^2$. And it was proven that $LE_k(G)=\sum_{i=1}^nd_i(d_i+1)$ which implies that $LE_k(G)$ is only related to
vertex degrees of $G$. Moreover, Gutman and Zhou \cite{GZ} provided
another definition of Laplacian energy of $G$ as follows.
$$LE_g(G)=\sum_{i=1}^n|\mu_i-\frac{2m}{n}|.$$ This definition preserves
the main features of the original graph energy.

Since the energy of undirected graphs was generalized to the skew energy
of oriented graphs, the Laplacian energy also has a generalization in oriented graphs.
Let $G^\sigma$ be an oriented graph with skew adjacency matrix $S(G^\sigma)$.
Then $SL(G^\sigma)=D(G)-S(G^\sigma)$ is called the skew Laplacian matrix of $G^\sigma$.
Suppose that $\nu_1,\nu_2,\ldots,\nu_n$ are the eigenvalues of $SL(G^\sigma)$.
Analogously, there are two different definitions of skew Laplacian energy
basing on the eigenvalues of $SL(G^\sigma)$.

Similar to $LE_k(G)$, Adiga and Smitha \cite{AS} defined
the skew Laplacian energy of $G^\sigma$ as
$$SLE_k(G^\sigma)=\sum_{i=1}^n\nu_i^2.$$

The following theorem indicates that $SLE_k(G^\sigma)$ is also only
related to the vertex degrees of $G$, independent of orientation of $G^\sigma$.
\begin{thm}\cite{AS}
Let $G$ be an undirected graph with vertex degrees $d_1,d_2,\ldots,d_n$.
Then for any oriented graph $G^\sigma$ of $G$, we have
$$SLE_k(G^\sigma)=\sum_{i=1}^nd_i(d_i-1).$$
\end{thm}

Then we can immediately deduce the following corollaries.
\begin{cor}\cite{AS}
For any oriented graph, its skew Laplacian energy is an even integer.
\end{cor}
\begin{cor}\cite{AS}
Let $H$ be a proper subgraph of a connected graph $G$ with at least
three vertices. Let $H^\tau$ and $G^\sigma$ be any oriented graphs
of $H$ and $G$, respectively. Then $SLE_k(H^\tau)<SLE_k(G^\sigma)$.
\end{cor}

An upper bound and a lower bound of $SLE_k(G^\sigma)$ are given in
the following theorem, and the extremal graphs are characterized therein.
\begin{thm}\cite{AS}
Let $G$ be a connected graph with $n\geq2$ vertices and $G^\sigma$
be any oriented graph of $G$. Then we have
$$2n-4\leq SLE_k(G^\sigma)\leq n(n-1)(n-2).$$
Moreover, $SLE_k(G^\sigma)=n(n-1)(n-2)$ if and only if $G^\sigma$ be
any oriented graph of the complete graph $K_n$ and $SLE_k(G^\sigma)=2n-4$
if and only if $G^\sigma$ be any oriented graph of the path $P_n$.
\end{thm}

Note that the above definition of skew Laplacian energy of an oriented graph
is independent of its orientation, which does not reflect the adjacency of
this oriented graph. Being aware of this, later Adiga and Khoshbakht \cite{AK}
gave another definition as follows.
$$SLE_g(G^\sigma)=\sum_{i=1}^n|\nu_i-\frac{2m}{n}|.$$
They also established some bounds analogous to these of $LE_g(G)$.
\begin{thm}\cite{AK}
Let $G$ be a graph with $n$ vertices, $m$ edges and degree sequence
$d_1,d_2,\ldots,d_n$. Let $G^\sigma$ be an oriented graph of $G$ with
skew Laplacian matrix $SL(G^\sigma)$. Assume that $\nu_1,\nu_2,\ldots,\nu_n$
are the eigenvalues of $SL(G^\sigma)$. Let $\gamma_i=\nu_i-\frac{2m}{n}$
and $|\gamma_1|\leq|\gamma_2|\leq\ \cdots\leq|\gamma_n|=k$. Then
\begin{enumerate}[\indent(i)]
\item $2\sqrt{M}\leq SLE_g(G^\sigma)\leq\sqrt{2M_1n}$;
\item $SLE_g(G^\sigma)\leq k+\sqrt{(n-1)(2M_1-k^2)}$;
\item If $G^\sigma$ has no isolated vertices, then $SLE_g(G^\sigma)\leq2M_1$, where $M=-m+\frac{1}{2}\sum_{i=1}^n(d_i-\frac{2m}{n})^2$ and $M_1=M+2m=m+\frac{1}{2}\sum_{i=1}^n(d_i-\frac{2m}{n})^2$.
\end{enumerate}
\end{thm}

Moreover, Kissani and Mizoguchi \cite{KM} introduced a different
Laplacian energy for directed graphs in 2010, where only the
out-degrees of vertices are considered rather than both the
out-degrees and in-degrees. The definition allows the directed
graphs to contain loops and dicycles with length $2$, but
its does not make use of the in-adjacency information of a digraph.
We omit the details there. In the following we focus on a new definition
\cite{CLS} of skew Laplacian energy of oriented graphs, which makes up for this shortage.

Let $G$ be a simple undirected graph with vertex set $V(G)=\{v_1,v_2,\ldots,v_n\}$.
Let $A(G)$ be the adjacency matrix of $G$ and $D(G)$ be the diagonal matrix
of vertex degrees of $G$. Let $G^\sigma$ be an oriented graph of $G$ with
skew adjacency matrix $S(G^\sigma)$. Denote by $d_i^+$ and $d_i^-$
the out-degree and in-degree of the vertex $v_i$ in $G$, respectively.
Define $\widetilde{D}(G^\sigma)=\mbox{diag}(d_1^+-d_1^-,d_2^+-d_2^-,\ldots,d_n^+-d_n^-)$, Inspired by the definition of laplacian matrix of undirected graphs, Cai, Li and Song
\cite{CLS} gave a new definition of skew Laplacian matrix as follows.
$$\widetilde{SL}(G^\sigma)=\widetilde{D}(G^\sigma)-S(G^\sigma).$$

Let $\{\eta_1,\eta_2,\ldots,\eta_n\}$ be the spectrum of $\widetilde{SL}(G^\sigma)$.
Note that $\widetilde{SL}(G^\sigma)$ is not (skew) symmetric. Thus its eigenvalues
may be imaginary numbers. Then the skew Laplacian energy \cite{CLS} of the
oriented graph $G^\sigma$ is defined as
$$SLE(G^\sigma)=\sum_{i=1}^n |\eta_i|.$$

Note that if an oriented graph is Eulerian, i.e. for every vertex
$v_i$, $d_i^+=d_i^-$, then its skew Laplacian energy equals its skew energy.
\begin{thm}\cite{CLS}
If $G$ is an Eulerian oriented graph, then $SLE(G^\sigma)=E_S(G^\sigma)$.
\end{thm}

In the following, we consider the bounds of this skew Laplacian energy.
\begin{thm}\cite{CLS}
Let $G$ be a graph with $n$ vertices, $m$ edges and $p$ components.
Let $G^\sigma $ be an oriented graph of $G$. Assume that $d_i^+$ ($d_i^-$)
is the out-degree ( in-degree) of the vertex $v_i$ in $G^\sigma$. Then
$$2\sqrt{|M|}\leq SLE(G^\sigma)\leq \sqrt{2M_1(n-p)},$$
where $M=-m+\frac{1}{2}\sum_{i=1}^n(d_i^+-d_i^-)^2$ and
$M_1=M+2m=m+\frac{1}{2}\sum_{i=1}^n(d_i^+-d_i^-)^2$. Moreover, there bounds are sharp.
\end{thm}
\begin{cor}\cite{CLS}
Let $G$ be a graph with $p$ components $C_1,C_2,\ldots,C_p$ and $G^\sigma$
be an oriented graph of $G$. If $SLE(G^\sigma)= \sqrt{2M_1(n-p)}$,
then each component $C_i$ is Eulerian with odd number of vertices.
\end{cor}
\begin{cor}\cite{CLS}
$SLE(G^\sigma)\leq \sqrt{2M_1n}$.
\end{cor}
\begin{cor}\cite{CLS}
If $G^\sigma$ has no isolated vertices, then $SLE(G^\sigma)\leq 2M_1$.
\end{cor}

\end{document}